\documentclass[reqno]{amsart} 
\usepackage{amsfonts, amsmath, amsthm, amssymb, latexsym, graphicx, geometry, xcolor, colortbl,  mathtools}
\usepackage{color}

\def\cth{{\rm coth}\>}

\def\b{\mathbb }
\def\cal{\mathcal }

\theoremstyle{plain}
\newtheorem{theorem}{Theorem}[section]
\newtheorem{corollary}[theorem]{Corollary}
\newtheorem{lemma}[theorem]{Lemma}
\newtheorem{proposition}[theorem]{Proposition}
\newtheorem{def-theorem}[theorem]{Definition and Theorem}
\theoremstyle{definition}
\newtheorem{definition}[theorem]{Definition}

\newtheorem{remark}[theorem]{Remark}

\newtheorem{example}[theorem]{Example}

\numberwithin{equation}{section}

\title[Bessel  and Dunkl processes with drift]{Bessel  and Dunkl processes with drift}

\author{ Michael Voit} 
\address{Fakult\"at Mathematik, Technische Universit\"at Dortmund,
          Vogelpothsweg 87,
          D-44221 Dortmund, Germany}
\email{ michael.voit@math.tu-dortmund.de}

\subjclass[2010]{Primary 60B20; Secondary 60F05, 60F15,  33C45, 60J60, 60K35,  70F10, 82C22}
\keywords{Multivariate Bessel functions, Dunkl processes, Bessel processes with drift, Dyson Brownian motions with drift, modified moments,
   strong laws of large numbers, central limit theorems.}

\begin{document}
\date{\today}

\begin{abstract} For certain discrete multiplicity parameters $k\ge0$, multivariate (Dunkl-)Bessel processes on Weyl chambers $C$
  associated with root systems appear as projections of Brownian motions without drift on Euclidean spaces $V$, and
  the associated transition densities  can be described in terms of multivariate Bessel functions;
  the most prominent examples are Dyson Brownian motions.
  The projections of Brownian motions on $V$ with drift are also Feller diffusions on $C$, and their
  transition densities and  their generators can be again described via these  Bessel functions. These processes are called
  Bessel processes with drift. In this paper we  construct these Bessel processes   with drift for arbitrary parameters $k\ge 0$. Moreover,
  this construction works also  for Dunkl processes.
  We study some features of these processes with drift like their radial parts, a Girsanov theorem, moments and associated martingales,
  strong laws of large numbers,
  and central limit theorems.
\end{abstract}

\maketitle

\section{Introduction}

It is well-known (see e.g. \cite{Ki, PY, RY}) that for a Brownian motion $(B_t)_{t\ge0}$ on $\mathbb R^n$ with start in $0\in\mathbb R^n$, its radial part $(\|B_t\|)_{t\ge0}$ (with the usual Euclidean norm on $\mathbb R^n$)
is a Bessel process with
dimension parameter $n$, where the transition probabilities can be written down by using spherical Bessel functions,  where 
these spherical Bessel functions and Bessel processes exist for all   parameters $n\in [1,\infty[$.
Moreover, it is also known that for arbitrary drift vectors $\tilde\lambda\in\mathbb R^n$, the radial parts $(\|B_t+t\tilde\lambda \|)_{t\ge0}$
of the Brownian motions $(B_t+t\tilde\lambda )_{t\ge0}$ with drift $\tilde\lambda$ are also diffusions whose distributions are the square-roots of non-central 
$\chi^2$-distributions.  Again, these diffusions can be defined for all   $n\in [1,\infty[$ in terms of the associated
stochastic differential equations, the associated generators, or by explicit transition probabilities by using the spherical Bessel functions.
We shall call these processes one-dimensional Bessel processes with drift on $[0,\infty[$.
    More precisely, as a special case of Subsection \ref{geometric-b}
    below for $N=1$, and with  $k:=(n-1)/2\ge0$, $\lambda:=\|\tilde\lambda\|\ge0$, we consider the
 spherical Bessel functions
 \begin{equation}\label{def-1dim-besselfkt-intro} j_{k-1/2}(z):=
   \Gamma(k+1/2)\sum_{l=0}^\infty \frac{(-1)^l (z/2)^{2l}}{l! \cdot \Gamma(l+k+1/2)} \quad\quad(z\in\mathbb C),
 \end{equation}
 and the positive functions $\tilde j_{k-1/2}(x):=j_{k-1/2}(ix)$ ($x\in\mathbb R$). Then
the Bessel processes with drift have the generators
\begin{equation}\label{Bessel-laplace-drift-gen-1dim}
  L_k^\lambda f(x):=
  \frac{1}{2} f^{\prime\prime}(x)  +
  \Bigl(\frac{k}{x}+\frac{\lambda \tilde j_{k-1/2}^\prime(\lambda x)}{ \tilde j_{k-1/2}(\lambda x)}\Bigr) f^{\prime}(x)
\end{equation}
for even   $f\in C^{(2)}(\mathbb R)$ with a reflecting boundary at $0$.
  These diffusions have the transition densities
   \begin{equation}\label{density-transition-bessel-drift-1dim}
  K_t^{\lambda,k}(x,A)=\frac{1}{2^{k-1/2}\Gamma(k+1/2)}\cdot  \frac{e^{-\lambda^2 t/2}}{t^{k+1/2}} \int_A e^{-(x^2+y^2)/(2t)}
 \frac{\tilde j_{k-1/2}(xy/t) \cdot\tilde j_{k-1/2}(y\lambda)}{\tilde j_{k-1/2}(x\lambda)}
\cdot y^{2k}\> dy
   \end{equation}
   for $t>0$, $x\in[0,\infty[$, and Borel sets $A\subset  [0,\infty[$.
   Hence the  associated transition operators appear,  up some normalizations,  as conjugations of the
transition operators without drift
   by multiplication operators which belong to the eigenfunctions  $y\mapsto j_{k-1/2}(y\lambda)$ of the
   generators $ L_k^0$ of the processes without drift.
   For the ``geometric cases'' $k:=(n-1)/2\ge0$ with  $n\in\mathbb N$, this  follows easily
   from the fact that these processes are projections of Brownian motions on  $\mathbb R^n$ with and without drift where this connection also holds; see e.g.~\cite{W, RP, DDMY, AuV}.
   This property is also preserved when we project  Brownian motions on  $\mathbb R^n$ with and without drift to other  orbit spaces
   which appear when
   suitable compact subgroups of the orthogonal group $O(n)$ act on  $\mathbb R^n$.  Typical examples  appear for instance
for  $\mathbb F=\mathbb R, \mathbb C, \mathbb H$, where
the unitary group $U(N,\mathbb F)$ acts by conjugation on the vector space $H_N(\mathbb F)$ of $N\times N$ Hermitian matrices over $\mathbb  F$ where
the space of orbits can be identified with the space of  ordered spectra of these matrices, i.e., with the closed Weyl chamber 
$$C_N^A= \{x\in \mathbb R^N: \quad x_1\ge x_2\ge\ldots\ge x_N\}.$$ In this case it is well known from random matrix theory (see e.g.~\cite{AGZ, D, M, F})
that the projections of  Brownian motions on $H_N(\mathbb F)$ without drift are examples of Dyson Brownian motions.
Motivated from harmonic analysis, these processes are also known as (Dunkl)-Bessel processes on $C_N^A$
associated with the root system $A_{N-1}$ (see e.g.~\cite{RV2, CGY, GY1, GY2, GY3}). In fact, it is possible to define  Bessel processes on closed
Weyl chambers without drift
for
all root systems where the number of parameters, called multiplicities, depends on the root system. Like in the $A_{N-1}$-case, these processes
are ``geometric'', i.e., projection of  Brownian motions without drift on suitable Euclidean spaces, for some discrete  multiplicities.
Furthermore, motivated by the theory of Dunkl operators (see e.g.~\cite{D1, D2, DX, A, Ka, R1, R3a}), these Bessel processes on Weyl chambers $C$ can be
extended to Feller processes on $\mathbb R^N$ where these processes behave like  Bessel processes with a.s.~continuous paths, and where jumps between the Weyl chambers are added by
random reflections associated with the  root system. These processes are called Dunkl processes; see \cite{R2, RV1, RV2, CGY, GY1, GY2, GY3, De}.
These Dunkl processes on $\mathbb R^N$ then are related to the corresponding Bessel processes on $C$ by the fact, that $C$ may be seen
as the set of orbits,
when the finite Weyl group generated by the reflections of the root system acts on  $\mathbb R^N$, and that  the corresponding Bessel processes are just
projections of the Dunkl processes. These Dunkl processes have much in common with  Brownian motions on  $\mathbb R^N$ (in fact, Brownian motions appear
for the multiplicity $k=0$). For instance, they are
martingales on $\mathbb R^N$, they have the same time inversion property like Brownian motions,
there are an associated heat polynomials, and there is an associated Wiener chaos and so on; see  \cite{ RV1, RV2, CGY, GY1, GY2, GY3, De}.

The aim of this paper is to extend the known theory of multivariate Bessel and Dunkl processes (without drift)
in a systematic way to these processes with drift where the general definition of these processes is motivated by the geometric examples above.
We here mention that the definition of the Bessel processes will be completely canonical, while in the Dunkl case there are two canonical
constructions, namely by conjugations of the cases without drift by multiplication operators associated with either  the associated
Dunkl kernels or Bessel functions. We shall see that the conjugation by Dunkl kernels leads to processes, which have much nicer properties.
For these reasons we shall call these processes Dunkl processes with drift, while the processes in the second case,
which are quite close to Bessel processes with drift, will be called hybrid Dunkl-Bessel processes with drift.

This paper is organized as follows:
In Section 2 we recapitulate several known facts from Dunkl theory like Dunkl Laplacians,
the Dunkl kernels, Bessel functions associated with root systems, as well
as Dunkl and Bessel processes without drift. In the central Section 3 we then discuss several classes of examples how Brownian motions with drift
on  finite-dimensional vector spaces leads to diffusion processes on associated closed Weyl chambers which we then call Bessel processes with drift.
The explicit generators and transition probabilities (in terms of  Bessel functions associated with root systems) then motivate the definition of
general  Bessel processes with drift for arbitrary root systems and multiplicities $k\ge 0$.
Moreover, we also construct
two kinds of Dunkl processes with drift. In the first construction we use Dunkl kernels in the conjugations
which leads to processes which have a behavior which is quite different from that of the associated  Bessel processes with drift.
This is quite different from the case of these processes without drift. For the second construction we use Bessel functions in the conjugations.
These processes on $\mathbb R^N$ then will be much closer to the  Bessel processes with drift. In particular they will have the same jump behavior like 
classical Dunkl processes without drift. We shall call these processes hybrid Dunkl-Bessel processes with drift.
After the construction of  Bessel processes, Dunkl processes, and hybrid Dunkl-Bessel processes with drift we
show that the radial parts of these processes
are just one-dimensional Bessel processes with drift analogous to the known case of these processes without drift.
Section 4 is then devoted to some simple Girsanov-type result  which connects these processes with drift to processes without drift by a
suitable simple transformation of the measure. In Section 5 we then introduce modified moments for Dunkl processes with drift
which are adapted to  these processes from an algebraic point of view. These moments then will lead to a strong law of large numbers (SLLN) and a central limit theorem (CLT).
Unfortunately, at some stage there appear serious problems with respect to the asymptotic behavior of Dunkl kernels.
For this reason, some results will remain in some implicit form, and completely satisfying results will be given for the case $N=1$ only.
Finally, in Section 6 we study the SLLN and CLT for Bessel processes with drift.
For the  hybrid Dunkl-Bessel processes with drift there do not exist interesting limit results which are not covered by the Bessel cases.

\section{Basics on Dunkl theory}

In this section we  recapitulate the well-known
definitions of Dunkl operators, Dunkl kernels, Bessel functions, and Bessel   processes (without drift) on closed
Weyl chambers $C\subset \mathbb R^N$ as well of Dunkl processes on $\mathbb R^N$; see e.g. \cite{R1, R2, R3a, RV1, RV3, CGY, GY1, GY2, GY3,
  Ka, A, DX, D2, De}
for the background.

\subsection{Root systems and Dunkl operators}\label{subsection-bessel-processes}
Let $\mathbb R^N$ be equipped with  the usual scalar product   $\langle\,.\,,\,.\,\rangle$, the associated 2-norm $\|.\|$,
and the standard basis $e_1,\ldots,e_N$.
Let $R\subset \mathbb R^N\setminus\{0\}$ be a  root system on $\mathbb R^N$
i.e., $R$ is a finite set of vectors $\alpha$, such that the reflections
 $\sigma_\alpha$  on the  hyperplanes perpendicular to $\alpha$  
 generate a finite group of orthogonal transformations on  $\mathbb R^N$, where for each  $\alpha\in R$, $\sigma_\alpha(R)=R$ holds.
 The group  $W$ generated by the  $\sigma_\alpha$ is called the finite reflection group or  Weyl group  associated with $R$.
 We  mainly consider the following  root systems:
 \begin{align}\label{ex-rootsystems}
   A_{N-1} &= \{  \pm(e_i - e_j); \>\> 1\leq i < j \leq N\} \quad\quad (N\ge 2),\notag\\
  B_N &=\{\pm e_i,  \pm(e_i \pm e_j); \>\> 1\leq i,j \leq N,\>\> i < j  \} \quad\quad (N\ge 1),\\
D_N&= \{\pm (e_i\pm e_j): \quad 1\le i<j\le N\} \quad\quad (N\ge 2).
  \notag
   \end{align}
 For $R=A_{N-1}$, the Weyl group $W$ is the symmetric group ${\cal S}_N$ permuting all coordinates, and 
 in the  $B_N$-case,  $W$ is  the hyperoctahedral group, where in addition sign changes appear in all coordinates.
 
Now fix some positive subsystem $R_+\subset R$ where we only take the roots in  (\ref{ex-rootsystems})
 where   instead of the first $\pm$-signs only  the $+$-signs appear. Moreover, fix some multiplicity function $k:R\to[0,\infty[$, i.e.,
  a function which is invariant under the action of $W$ on $R$. For the root systems
  $A_{N-1}$ and $D_N$, this $k$ is  a single parameter $k\in [0,\infty[$, and, for
 $B_N$,   $k$  consists of 2 parameters 
  $(k_1,k_2)$  where  $k_1,$ and $k_2$ are the values on the roots $\pm e_i$ and $ e_i \pm e_j$ respectively.
      We  fix the closed Weyl chamber
      $$C:=\{ x\in \mathbb R^N: \quad \langle x, \alpha\rangle\ge 0 \quad\text{for all}\quad \alpha\in R_+\}$$
associated with  $R_+$. This chamber $C$     
consists of representatives of the orbits under the action of $W$ on  $\mathbb R^N$. In our examples, we 
have
\begin{align}\label{chambers}
  C_N^A&= \{x\in \mathbb R^N: \quad x_1\ge x_2\ge\ldots\ge x_N\},\\
 C_N^B &=\{x\in \mathbb R^N: \quad x_1\ge x_2\ge\ldots\ge x_N\ge 0 \},\notag\\
C_N^D&=\{x\in\mathbb R^N: \quad x_1\ge \ldots\ge x_{N-1}\ge |x_N|\}.
  \notag
\end{align}
We now fix $R$, $R_+$, $W$, $k$, and  $ C\subset \mathbb R^N$ and consider 
the associated Dunkl operators 
\begin{equation}\label{def-dunkl-op}
  T_\xi(k) f(x)   = \partial_{\xi}  f(x)+\sum_{\alpha\in R_+}
  k(\alpha)\cdot \langle\alpha,\xi\rangle \frac{1}{\langle\alpha, x\rangle} (f(x)-f(\sigma_\alpha(x)))
\end{equation} 
for  $\xi \in \mathbb R^N$, the directional derivatives $\partial_{\xi}$, and $f\in C^{(1)}(\mathbb R^N)$ where the $T_\xi(k) f$
are continuous on $\mathbb R^N$.
The  $T_\xi(k)$ form a  commuting family of  operators which are homogeneous of degree $-1$ on 
the space 
of all polynomials in $N$ variables such that for any  $r$-times continuously differentiable functions  $f$ on $\mathbb R^N$,  the
$T_\xi(k) f$ are $(r-1)$-times differentiable.

Moreover, 
for each  $w\in \b C^N$, the joint eigenvalue problem 
\begin{equation}\label{joint-Dunkl-eigenvalues}
  T_\xi(k) f\,=\, \langle\xi,w\rangle f \quad (\xi\in \b R^N) \quad\text{with}\quad  \quad f(0)=1
  \end{equation}
has a unique holomorphic solution $f(z) = E_k(z,w)$, called the 
Dunkl kernel. We list some well-known properties of $E_k$:

\begin{proposition}\label{properties-dunkl-kernel}
 For $x,y\in \mathbb C^N, \,\lambda\in \mathbb C$ and $g\in W$:
\begin{enumerate}\itemsep=-1pt
\item[\rm{(1)}] $E_k(x,y) = E_k(y,x).$ 
\item[\rm{(2)}]  $E_k(\lambda x,y) = E_k(x,\lambda y)$. 
\item[\rm{(3)}] $\overline{E_k(x,y)} = E_k(\overline x,\overline y).$
\item[\rm{(4)}]  $0<E_k( x,y) \le \max_{g\in W} e^{\langle x, gy\rangle}$ and  $|E_k( x,iy)|\le 1$ for $x,y\in \mathbb R^N$.
 \item[\rm{(5)}]  $E_k(gx,gy) = E_k(x,y)$.
\end{enumerate}
\end{proposition}

For $k\ge0$, the Dunkl kernels
$E_k$ have the following positive integral representation in terms of exponential functions:

\begin{theorem}\label{general-positive-integralrep-dunkl}
  For each $x\in\mathbb R^N$ there exists a unique probability measure $\mu_x\in M^1(\mathbb R^N)$ with
  $$E_k(x,\lambda)= \int_{\mathbb R^N} e^{\langle y,\lambda\rangle}\> d\mu_x(y) \quad\text{ for all}\quad \lambda\in \mathbb C^N.$$
  The measures $\mu_x$ have the following properties:
  \begin{enumerate}\itemsep=-1pt
\item[\rm{(1)}] The support $supp\>\mu_x $ is contained in the convex hull of the Weyl group orbit $\{gx: \> g\in W\}$ of $x$.
 \item[\rm{(2)}]  $x\in supp\>\mu_x $, and for $k>0$,  $\{gx: \> g\in W\}\subset  supp\>\mu_x $.
\item[\rm{(3)}] For $t>0$ and $x\in\mathbb R^N$,
  $\mu_{tx}$ is the image of $\mu_x$ under the mapping $y\mapsto ty$ on $\mathbb R^N$.
  \end{enumerate}
\end{theorem}

\begin{proof} The existence of the  $\mu_x$ is shown in \cite{R3}. For the first statement in (2) see \cite{dJ},
  and for the second one in (2) see \cite{GR}. Finally,
    (3)  follows from
  $$\int_{\mathbb R^N} e^{\langle y,\lambda\rangle}\> d\mu_{tx}(y)=E_k(tx,\lambda)=E_k(x,t\lambda)=
  \int_{\mathbb R^N}  e^{\langle y,t\lambda\rangle}\> d\mu_x(y) =
\int_{\mathbb R^N}  e^{\langle ty,\lambda\rangle}\> d\mu_x(y) .$$
  \end{proof}

  In general, no explicit descriptions are known about the measures $\mu_x$  except for the trivial case $k=0$ and
  for the case  $N=1$ from \cite{R1}:

\begin{proposition}\label{int-rep-Dunkl-1dim-Prp}
For $N=1$, $k>0$, and $\lambda,x\in\mathbb R$,
\begin{equation}\label{int-rep-Dunkl-1dim}
E_k(x,\lambda)=\frac{\Gamma(k+1/2)}{\Gamma(1/2)\Gamma(k)}\int_{-1}^1 e^{x\lambda t} (1-t)^{k-1}(1+t)^k \> dt.
  \end{equation}
\end{proposition}

  Theorem \ref{general-positive-integralrep-dunkl}(2) implies the following lower bounds for $E_k$:

  \begin{corollary}\label{est-below} For each $\lambda\in\mathbb R^N$ there are constants $c_1>0$ and $c\in\mathbb R^N$
    such that $E_k(x,\lambda)\ge c_1 e^{\langle c,x\rangle}$ for $x\in\mathbb R^N$.
\end{corollary}

We next define the  Bessel functions
$$ J_k(z,w) := \frac{1}{|W|}\sum_{g\in W} E_k(z,gw)$$
as symmetrizations of the $E_k$. Then all assertions of Proposition \ref{properties-dunkl-kernel},
Theorem \ref{general-positive-integralrep-dunkl}, and Corollary  \ref{est-below} hold also for $J_k$,
where instead of \ref{properties-dunkl-kernel}(5),
$J_k$ even is  $W$-invariant in both arguments separately. Moreover,
 for  $w\in \mathbb C^N$, $g(z):= J_k(z,w)$ is the unique holomorphic solution of the "Bessel system''
\begin{equation}\label{besselsystem}
p(T_{e_1}(k),\ldots,T_{e_N}(k)) g\, = \,p(w)g \quad\forall\, p\in \mathcal P^W \quad\text{with}\quad g(0)=1
\end{equation}
with the algebra
$\mathcal P^W$  of $W$-invariant polynomials in $N$ variables.
Restricted to sufficiently smooth, $W$-invariant functions on $\mathbb R^N$,
the operators $p(T_{e_1}(k),\ldots,T_{e_N}(k))$  are differential operators.
In particular, for $p(x)=x_1^2+\ldots+x_N^2$, one obtains the Dunkl-Laplacian $\Delta_k$ which has the form
\begin{equation}\label{Dunkl-laplace}
 \Delta_k f(x) = \Delta f(x)  + \,2\sum_{\alpha\in R_+}
k(\alpha)\Bigl(
\frac{\langle\nabla f(x),\alpha\rangle}{\langle\alpha ,x\rangle} -\frac{ \|\alpha\|^2}{2}
\frac{f(x)-f(\sigma_\alpha x)}{\langle\alpha , x\rangle^2}\Bigr).
\end{equation}
Note that for $W$-invariant functions $f$ the reflection parts at the end of \eqref{Dunkl-laplace} disappear.

We  need the following product rule which follows  either by  elementary calculus or
from Lemma 3.1 in \cite{Ve} by polarization and the fact that in \cite{Ve}
all roots have the normalization $\|\alpha\|^2=2$.

\begin{lemma}\label{prod-rule}
  For all functions $f,g\in C^{(2)}(\mathbb R^N)$,
  \begin{align}
    \Delta_k (fg)(x)= f(x) \Delta_kg(x)+ & g(x) \Delta_kf(x)+ 2\langle \nabla f(x), \nabla g(x)\rangle \notag\\&+
    \sum_{\alpha\in R_+}k(\alpha)\|\alpha\|^2 \frac{(f(x)-f(\sigma_\alpha x))(g(x)-g(\sigma_\alpha x))}{\langle\alpha , x\rangle^2}.\notag\end{align}
  Moreover, if $f$ or  $g$ is $W$-invariant, then the reflection parts on the right hand side disappear.
\end{lemma}

It is well-known (see \cite{R2, RV1, RV2}) that  $\frac{1}{2}\Delta_k$ is the generator of a Feller process
$(X_{t,k})_{t\ge0}$ on $\mathbb R^N$ with RCLL paths; the factor $1/2$ here is used in order to regard these processes as generalizations of Brownian motions with the usual time normalization in probability. These processes have
 the transition  probabilities
\begin{equation}\label{density-transition-dunkl-gen}
K_t^k(x,A)=c_k \int_A \frac{1}{t^{\gamma+N/2}} e^{-(\|x\|^2+\|y\|^2)/(2t)} E_k(\frac{x}{\sqrt{t}}, \frac{y}{\sqrt{t}}) 
\cdot w_k(y)\> dy
\end{equation}
for $t>0$, $x\in \mathbb R^N$, and Borel sets $A\subset \mathbb R^N$
with  the  $W$-invariant  weight function
\begin{equation}\label{def-wk-A-general}
   w_k(x) := \prod_{\alpha\in R_+} |\langle\alpha,x\rangle|^{2k(\alpha)}.
\end{equation}
 which is homogeneous of  degree $2\gamma$ with 
$\gamma:=\sum_{\alpha\in R_+}k(\alpha)$, and where
 $c_k>0$ is a known  normalization depending on $k$ and $R$ only. For instance, by \cite{M},
 for the root systems $A_{N-1}$ and $B_N$ and the multiplicities $k$ and $(k_1,k_2)$, these constants are
 \begin{align}\label{ck-a-b}
    c_k^A=&\frac{1}{(2\pi)^{N/2}} \cdot\prod_{j=1}^{N}\frac{\Gamma(1+k)}{\Gamma(1+jk)}, \\
 c_k^B=&\frac{1}{2^{N(k_1+(N-1)k_2+1/2)}} \cdot\prod_{j=1}^{N}\frac{\Gamma(1+k_2)}{\Gamma(1+jk_2)\Gamma(\frac{1}{2}+k_1+(j-1)k_2)}
\end{align}
respectively.
The processes $(X_{t,k})_{t\ge0}$
 are called  Dunkl processes associated with $R$ and $k$.

 Furthermore, if we identify $C$ with the space of all $W$-orbits in $ \mathbb R^N$, then the operator
 \begin{equation}\label{Bessel-laplace}
  L_k f(x):= \frac{1}{2}\Delta_k f(x)=\frac{1}{2}\Delta f(x)  + \sum_{\alpha\in R_+}
k(\alpha)
\frac{\langle\nabla f(x),\alpha\rangle}{\langle\alpha ,x\rangle}
\end{equation}
for $W$-invariant $f\in C^{(2)}(\mathbb R^N)$ is the generator of a Feller diffusion
$(X_{t,k})_{t\ge0}$ on $C$ with continuous paths with reflecting boundaries such that  the transition  probabilities satisfy
\begin{equation}\label{density-transition-bessel-gen}
K_t^k(x,A)=|W|c_k \int_A \frac{1}{t^{\gamma+N/2}} e^{-(\|x\|^2+\|y\|^2)/(2t)} J_k(\frac{x}{\sqrt{t}}, \frac{y}{\sqrt{t}}) 
\cdot w_k(y)\> dy
\end{equation}
for $t>0$, $x\in C$, and Borel sets $A\subset C$. We call these processes Bessel processes on $C$ associated with $R$ and $k$.
In some references these processes are called 
 Dunkl-Bessel  or invariant Dunkl processes.

 For the root systems $A_{N-1}, B_N, D_N$ the generators  $L_k$ in \eqref{Bessel-laplace} and the weights $w_k$ are as follows:
 In the $A$-case we have
 \begin{equation}\label{Bessel-laplace-a}
  L_k^A f(x):= \frac{1}{2}\Delta f(x)  + \sum_{1\le i<j\le N}
k
\frac{ f_{x_i}(x)-f_{x_j}(x)}{x_i-x_j},   \quad  w_k^A(x)=\prod_{1\le i<j\le N}(x_i-x_j)^{2k}
\end{equation}
with $k\ge0$ and $\gamma_A=kN(N-1)/2$, i.e., Bessel processes of type $A_{N-1}$ are Dyson Brownian motions as studied e.g.~in \cite{AGZ}.
 In the $B$-case we have for $k_1,k_2\ge0$
 \begin{equation}\label{Bessel-laplace-b}
   L_k^B f(x):= \frac{1}{2} \Delta f(x)  +
 k_2 \sum_{i=1}^N \sum_{j\ne i} \Bigl( \frac{1}{x_i-x_j}+\frac{1}{x_i+x_j}  \Bigr)
 \frac{\partial}{\partial x_i}f (x)
 \quad + k_1\sum_{i=1}^N\frac{1}{x_i}\frac{\partial}{\partial x_i}f(x)\end{equation}
 and
 \begin{equation}\label{def-wk-b}
   w_k^B(x):= \prod_{1\le i<j\le N}(x_i^2-x_j^2)^{2k_2}\cdot \prod_{i=1}^N x_i^{2k_1}, \quad  \gamma_B=k_2N(N-1)+k_1N.\end{equation}
 Finally, in the $D$-case, for $k\ge0$,
\begin{equation}\label{Bessel-laplace-d} L_k^Df(x):= \frac{1}{2} \Delta f(x) +
 k \sum_{i=1}^N \sum_{j\ne i} \Bigl( \frac{1}{x_i-x_j}+\frac{1}{x_i+x_j}  \Bigr)
 \frac{\partial}{\partial x_i}f(x) \end{equation}
and
\begin{equation}\label{def-wk-d}w_k^D(x):= \prod_{i<j}(x_i^2-x_j^2)^{2k}, \quad  \gamma_D:= kN(N-1).  \end{equation}

\section{Bessel and Dunkl processes with drift: Motivation and definition}

In this section we  introduce general multivariate Bessel and Dunkl processes with drift.
In the Bessel case, these processes are motivated by the fact that for particular ``geometric'' multiplicities $k\ge0$
these processes
appear as projections of Brownian motions with drift on  associated Euclidean spaces $V$.
For the case of Bessel processes without drift, this is well known; see e.g. \cite{Ki, RY, BH} and references there
for the case $N=1$, \cite{AGZ} in the $A_{N-1}$-case of Dyson Brownian motions, and
\cite{R4, RV2} and references there for further examples.  It was shown in \cite{AuV} (see also \cite{RV2})
that in these geometric situations, projections from $V$ onto the corresponding Weyl chambers $C$ of  Brownian motions with
arbitrary constant drift 
vectors are also Feller diffusions, where the associated transition densities and their generators can be described
in terms of the associated data of the original Bessel processes (without drift) and the associated Bessel functions.
Motivated by this connection, we call these Feller diffusions Bessel processes with drift on $C$.
This connection between Bessel processes with and without drift works also for arbitrary Bessel processes associated with root systems
and multiplicities $k\ge0$. We show that such a connection is also available for Dunkl processes on $ \mathbb R^N$.

\subsection{Bessel processes in the geometric setting}
Let 
$(V, \langle\,.\,,\,.\,\rangle)$ be a 
Euclidean vector space of finite dimension $n$ and
$K\subset O(V)$ a  compact subgroup  of the orthogonal group of $V$.
 Consider the exponential functions
 $$e_\lambda(x):=e^{\langle \lambda,x\rangle}  \quad\quad(x,\lambda\in V)$$
 and  the space  $V^K:=\{K.x: \> x\in V\}$  of all $K$-orbits $K.x:=\{g(x): \> g\in K\}$ in $V$.
 The space $V^K$ is a locally compact Hausdorff space w.r.t.~the quotient topology.
Let $p:V\to V^K$ be the canonical  continuous and open projection,
and  denote the associated push forward of measures also by $p$. We then have the following  observation;
see \cite{RV2, AuV}:

\begin{lemma}\label{modified-projection1} Let  $\lambda\in V$. Then, for $x\in V$, and the normalized Haar measure $dk$ of $K$,
  $$\alpha_\lambda(p(x)):=\int_K e_\lambda(k(x)) \> dk$$
  is a well-defined continuous function on $V^K$, and for all bounded Borel measures
  $$\mu\in M_{e,K}(V):=\{\mu\in M_b(V):\> \mu \quad K\text{-invariant, and }\quad e_{\tilde\lambda}\cdot \mu
  \quad \text{is bounded for all}\quad \tilde\lambda\in V\}$$
    we have
  \begin{equation}\label{general-image}
    p(e_\lambda \cdot \mu) = \alpha_\lambda \cdot p(\mu).\end{equation}
\end{lemma}

Moreover, by \cite{AuV}, we have:

\begin{proposition}\label{modified-projection}
  Let $\lambda\in V$ be some drift vector, and $(B_t)_{t\ge0}$ a Brownian motion (without drift) on $V$ starting in $0\in V$ where, by definition,
  the distributions $d\mu_t(x)= \frac{1}{(2\pi t)^{n/2}} e^{-\|x\|^2/(2t)}dx$ of $B_t$ are $K$-invariant by construction.
  Then $(B_t+t\lambda)_{t\ge0}$ is a Brownian motion  on $V$ with drift $\lambda$ where, clearly,  $B_t+t\lambda$ has the distribution
  $\mu_t^\lambda:= e^{-\|\lambda\|^2 t/2} e_\lambda \> \mu_t$ for $t\ge0$. The projection $(p(B_t+t\lambda))_{t\ge0}$ of this Brownian motion
  with drift is a Feller diffusion on $V^K$. Moreover:
  \begin{enumerate}
  \item[\rm{(1)}] The transition probabilities $p_t^\lambda(x,dy)$ and  $p_t^0(x,dy)$ of the projections of $(B_t+t\lambda)_{t\ge0}$ and of
    $(B_t=B_t+t\cdot 0)_{t\ge0}$ respectively with $t\ge 0$ and $x\in V^K$ are related by
    $$ dp_t^\lambda(x,dy)= \frac{e^{-\|\lambda\|^2 t/2}}{\alpha_\lambda(x) }\cdot \alpha_\lambda(y) dp_t^0(x,dy) \quad\quad(x\in V^K).$$
     \item[\rm{(2)}]  The generators $L^\lambda$ and $L^0$ of the transition semigroups of the projections of $(B_t+t\lambda)_{t\ge0}$ and of
  $(B_t=B_t+t\cdot 0)_{t\ge0}$ respectively are related by
\begin{equation}\label{trafo-generators}
  L^{\lambda}f(y)=  - \frac{\|\lambda\|^2}{2} f(y) +\frac{1}{\alpha_\lambda(y)} L^0(\alpha_\lambda f)(y)
\end{equation}
for all 
$f\in C_0(V^K)$ which are contained in  the domain $D(L^\lambda)$ of $L^\lambda$.
\end{enumerate}
\end{proposition}

We now apply this result to pairs $(V,K)$ where $V^K$ can be identified with a closed Weyl chamber
$C\subset\mathbb R^N$, and where the functions $\alpha_\lambda$ on $C$ are Bessel
functions associated with an associated root system $R\subset \mathbb R^N$. We discuss the examples for  $R=A_{N-1}, B_N,D_N$ separately:

\subsection{Geometric cases of type A}\label{geometric-a}
Take one of the fields $\mathbb F=\mathbb R, \mathbb C$ or the skew-field $\mathbb F=\mathbb H$ of quaternions
with the real dimension $d=1,2,4$ respectively. Let $V$ be the vector space $H_N(\mathbb F)$ (over $\mathbb R$) of all Hermitian $N\times N$-matrices
with the trace scalar product $\langle A,B\rangle:=\Re tr(A\bar B^{T})$. Then $K:=O(N,\mathbb F)$
acts on $V$ as a compact group of orthogonal transformations on  $H_N(\mathbb F)$ by conjugations, and two matrices  in
$H_N(\mathbb F)$ are in the same orbit precisely if they have the same real eigenvalues with the same multiplicities.
We thus can identify $V^K$ with the Weyl chamber $C_N^A$ such that $(x_1,\ldots,x_N)\in C_N^A$ corresponds to the orbit of the matrix
$diag(x_1,\ldots,x_N)\in  H_N(\mathbb F)$. Denote the associated projection  as above by $p:V\to C_N^A$ and put $k:=d/2=1/2,1,2$.
Then, for each ``drift vector'' $\Lambda\in H_N(\mathbb F)$, the function $\alpha_\Lambda$ on $ C_N^A$ from Lemma \ref{modified-projection1}
is just the Bessel function $J_k^A(.,p(\Lambda))$. Moreover, for a Brownian motion $(B_t)_{t\ge 0}$ on $H_N(\mathbb F)$, its projection 
$(p(B_t))_{t\ge 0}$ on $C_N^A$ is a Bessel process of type $A_{N-1}$ with this parameter $k$ starting in $0\in C_N^A$  with generator
\eqref{Bessel-laplace-a} and transition probabilities \eqref{density-transition-bessel-gen}.
These processes are also known as Dyson Brownian motions; see e.g. \cite{AGZ} and references there.

Now fix a drift $\Lambda\in H_N(\mathbb F)$ and put $\lambda:=p(\Lambda)\in C_N^A$.
 Proposition \ref{modified-projection}(1) and
\eqref{density-transition-bessel-gen} show that then the
diffusion $(X_t^{\lambda,k}:= p(B_t+t\Lambda))_{t\ge0}$ on $C_N^A$ with start in 0 has the transition probabilities
\begin{equation}\label{density-transition-bessel-a-drift}
  K_t^{\lambda,k}(x,A)=\frac{N!\cdot c_k\cdot e^{-\|\lambda\|^2 t/2}}{t^{\gamma+N/2}} \int_A e^{-(\|x\|^2+\|y\|^2)/(2t)}
 \frac{ J_k(\frac{x}{\sqrt{t}}, \frac{y}{\sqrt{t}}) \cdot J_k(y,\lambda)}{J_k(x,\lambda)}
\cdot \prod_{1\le i<j\le N}(y_i-y_j)^{2k}\> dy
\end{equation}
for $t>0$, $x\in C_N^A$ and Borel sets $A\subset  C_N^A$. Furthermore,  Proposition \ref{modified-projection}(2), \eqref{Bessel-laplace-a},
 Lemma \ref{prod-rule}
 (with the $W$-invariant function $g(y):= J_k(y,\lambda)$), and
 $\Delta_kJ(.,\lambda) = \|\lambda\|^2J(.,\lambda) $ ensure that this diffusion has the generator
\begin{equation}\label{Bessel-laplace-a-drift}
  L_k^{A,\lambda} f(x)= \frac{1}{2}\Delta f(x)  + \sum_{1\le i<j\le N}
k\frac{ f_{x_i}(x)-f_{x_j}(x)}{x_i-x_j} + \frac{\langle \nabla_x J_k(x,\lambda), \nabla f(x)\rangle }{J_k(x,\lambda)}
\end{equation}
for $W$-invariant $f\in C^{(2)}(\mathbb R^N)$ where $\nabla_x$ denotes the gradient w.r.t.~$x$.

The result for the root system $B_N$ is similar:

\subsection{Geometric cases of type B}\label{geometric-b}
Fix one of the (skew-)fields $\mathbb F=\mathbb R, \mathbb C, \mathbb H$  with the real dimension $d=1,2,4$ as above.
For integers $M\ge N\ge 1$
 consider the vector spaces $V:=M_{M,N}(\mathbb F)$ of all $M\times N$-matrices over $\mathbb F$ with
 the real dimension $dMN$ and the scalar product $\langle A,B\rangle:=\Re tr(A\bar B^{T})$.

 For  $A\in M_{M,N}(\mathbb F)$ consider the matrix $A^*:=\bar A^T\in M_{N,M}(\mathbb F)$ as well as  the matrix 
 $A^*A$ in the cone $\Pi_N(\mathbb F)$
of all $N\times N$ positive semidefinite matrices over  $\mathbb F$. 
Furthermore, consider the spectral mapping $\sigma_N:\Pi_N(\mathbb F)\to C_N^B$ which assigns to each  matrix in $\Pi_N(\mathbb F)$
its  ordered spectrum.
We now consider the projection $$p: M_{M,N}(\mathbb F)\to C_N^B, \quad A\mapsto \sqrt{\sigma_N(A^*A)}$$
where the symbol $\sqrt{.}$ means taking square roots in all coordinates.

It is well known that for a  Brownian motion
 $(B_t^M)_{t\ge0}$ on  $M_{M,N}(\mathbb F)$ , the process  $((B_t^M)^*B_t^M)_{t\ge0}$ is a Wishart process on $\Pi_N(\mathbb F)$  with shape parameter $M$
(see \cite{Bru,  DDMY} and references there
for details on Wishart processes), and that $(p(B_t^M))_{t\ge0}$ is a Bessel process on  $C_N^B$ of type $B_N$ with multiplicities
$k=(k_1,k_2):=((M-N+1)\cdot d/2-1/2, d/2)$ (see e.g. \cite{R4, RV2}).

For Brownian motions with drift we now proceed as in the preceding example.
Fix a drift matrix  $\Lambda\in  M_{M,N}(\mathbb F)$ and put $\lambda:=p(\Lambda)\in C_N^B$.
By Proposition \ref{modified-projection}(1) and
\eqref{density-transition-bessel-gen},  the
diffusion $(X_t^{\lambda,k}:= p(B_t^M+t\Lambda))_{t\ge0}$ on $C_N^B$ then has the transition probabilities
\begin{equation}\label{density-transition-bessel-b-drift}
  K_t^{\lambda,k}(x,A)=\frac{2^N N!\cdot c_k\cdot e^{-\|\lambda\|^2 t/2}}{t^{\gamma+N/2}} \int_A e^{-(\|x\|^2+\|y\|^2)/(2t)}
 \frac{ J_k(\frac{x}{\sqrt{t}}, \frac{y}{\sqrt{t}}) \cdot J_k(y,\lambda)}{J_k(x,\lambda)}
\cdot w_k^B(y)\> dy
\end{equation}
for $t>0$, $x\in C_N^B$ and $A\subset  C_N^B$ with $ w_k^B$ as in \eqref{def-wk-b}.
Furthermore,  by the same arguments as before, this diffusion has the generator
\begin{align}\label{Bessel-laplace-b-drift}
  L_k^{B,\lambda} f(x)= \frac{1}{2}\Delta f(x)  &+
 k_2 \sum_{i=1}^N \sum_{j\ne i} \Bigl( \frac{1}{x_i-x_j}+\frac{1}{x_i+x_j}  \Bigr)
 \frac{\partial}{\partial x_i}f (x)\notag\\
 &+ k_1\sum_{i=1}^N\frac{1}{x_i}\frac{\partial}{\partial x_i}f(x)
 + \frac{\langle \nabla_x J_k(x,\lambda), \nabla f(x)\rangle }{J_k(x,\lambda)}
\end{align}
for $W$-invariant $f\in C^{(2)}(\mathbb R^N)$.

Besides these examples there is a further geometric example which leads to a Bessel process of type $B_N$:
 For $N\ge 2$, consider the vector space $V:=Skew(2N+1,\mathbb R)$ of all skew-symmetric real matrices of dimension $2N+1$ on which the compact group
$K:=SO(2N+1,\mathbb R)$ acts by conjugation. $V^K$ can be identified with $C_N^B$ in a canonical way;
 see e.g. \cite{T} and  \cite{M}
 for  the random matrix context. Projections of Brownian motions on  $V$ with and without drift then again lead to
 processes of type $B_N$ as above with $k_1=k_2=1$; see also Section 4.1 in \cite{AuV} for further details.

\subsection{Geometric cases of type D}\label{geometric-d} For the root system $D_N$ there also exists a similar example with $V:=Skew(2N,\mathbb R)$
 on which  $K=SO(2N)$ acts by conjugation. Here, $V^K$ can be identified with
 $C_N^D$, and projections of Brownian motions without drift on $V$ are  Bessel processes of type $D_N$ with $k=1$;
 see Appendix C in \cite{Kn} and for its connections to random matrices \cite{M}. The corresponding data for the
projections of Brownian motions with drifts have the same form as above; see  Section 4.2 in \cite{AuV}.

As all examples of Bessel processes and their modifications above have the same structure, these examples
 motivate the following
general construction of Bessel  processes with drift:

\begin{def-theorem}\label{def-bessel-p-drift}
  Let $R,R_+, W, C,$ and $ k\ge0$ be given as above, and
 let $\lambda\in C$ be a ``drift vector''. Then the operator
   \begin{equation}\label{Bessel-laplace-drift-gen}
  L_k^\lambda f(x):= \frac{1}{2}\Delta f(x)  + \sum_{\alpha\in R_+}
k(\alpha)
\frac{\langle\nabla f(x),\alpha\rangle}{\langle\alpha ,x\rangle}+\frac{\langle \nabla_x J_k(x,\lambda), \nabla f(x)\rangle }{J_k(x,\lambda)}
   \end{equation}
   for $W$-invariant $f\in C^{(2)}(\mathbb R^N)$ is the generator of a Feller diffusion on $C$ with reflecting boundaries.
   These diffusions have the transition densities
   \begin{equation}\label{density-transition-bessel-drift-gen}
  K_t^{\lambda,k}(x,A)=\frac{|W|\cdot c_k\cdot e^{-\|\lambda\|^2 t/2}}{t^{\gamma+N/2}} \int_A e^{-(\|x\|^2+\|y\|^2)/(2t)}
 \frac{ J_k(\frac{x}{\sqrt{t}}, \frac{y}{\sqrt{t}}) \cdot J_k(y,\lambda)}{J_k(x,\lambda)}
\cdot w_k(y)\> dy
\end{equation}
for $t>0$, $x\in C$ and  Borel sets $A\subset  C$ with the weight $ w_k$ in \eqref{def-wk-A-general}.
These diffusions  have modifications with continuous paths; they
will be called Bessel processes of type $R$ with multiplicity $k$ and drift $\lambda$.   
\end{def-theorem}

We  prove this result below together with a corresponding result for Dunkl processes with drift. In the case of Dunkl processes
however, there are two
``canonical'' definitions of such processes. In the first version we replace in all instances Bessel functions by Dunkl kernels.
We show below that these processes with drift  have less in common with the Bessel processes 
in Theorem \ref{def-bessel-p-drift} than in the case $\lambda=0$ without drift. However, it will turn out in  Section 5
that the name ``Dunkl processes with drift $\lambda$'' fits quite well to these processes. For this reason we study these processes 
more closely below. On the other hand, we can define Dunkl processes with drift by replacing the Bessel functions
by Dunkl kernels only partially. It will turn out that these processes are closely related with the Bessel processes with drift
in Theorem \ref{def-bessel-p-drift} (see Remark \ref{connection-dunkl-bessel} for details).
On the other hand, these Dunkl-type processes with drift depend
on the Weyl group  orbit of the drift vector $\lambda$ and not on $\lambda$ itself. This means that  the interpretation of
$\lambda$ as a drift on $\mathbb R^N$ is not really correct here. Nevertheless, we call these processes hybrid  Dunkl-Bessel processes 
 with drift $\lambda$. The details of both definitions are as follows:

 \begin{def-theorem}\label{def-dunkl-bessel-p-drift}
    Let $R,R_+, W, C,$ as well as $ k\ge0$  and a ``drift'' $\lambda\in \mathbb R^N$ as above.
  \begin{enumerate}\itemsep=-1pt
  \item[\rm{(1)}] The operator
   \begin{align}\label{Dunkl-laplace-drift-gen}
  L^{\lambda, Dunkl}_k f(x):= &\frac{1}{2}\Delta f(x)  + \sum_{\alpha\in R_+}
k(\alpha)
\frac{\langle\nabla f(x),\alpha\rangle}{\langle\alpha ,x\rangle} + 
\frac{\langle \nabla_x E_k(x,\lambda), \nabla f(x)\rangle }{E_k(x,\lambda)}\\
&-\sum_{\alpha\in R_+}\frac{k(\alpha)\|\alpha\|^2}{2}\cdot \frac{ E_k(\sigma_\alpha x,\lambda)}{ E_k(x,\lambda)}\cdot
\frac{f(x)-f(\sigma_\alpha x)}{\langle\alpha , x\rangle^2} \notag
   \end{align}
   for  $f\in C^{(2)}(\mathbb R^N)$ is the generator of a Feller process on $\mathbb R^N$.
   These processes have the transition densities
   \begin{equation}\label{density-transition-dunkl-drift-gen}
  K_t^{\lambda,k,Dunkl}(x,A)=\frac{ c_k\cdot e^{-\|\lambda\|^2 t/2}}{t^{\gamma+N/2}} \int_A e^{-(\|x\|^2+\|y\|^2)/(2t)}
 \frac{ E_k(\frac{x}{\sqrt{t}}, \frac{y}{\sqrt{t}}) \cdot E_k(y,\lambda)}{E_k(x,\lambda)}
\cdot w_k(y)\> dy
\end{equation}
for $t>0$, $x\in \mathbb R^N$ and Borel sets $A\subset \mathbb R^N$ with  $ w_k$  as in \eqref{def-wk-A-general}.
These Feller processes have modifications such that all paths are right continuous with left limits (RCLL).
Such processes will be called Dunkl processes of type $R$ with multiplicity $k$ and drift $\lambda$.
\item[\rm{(2)}] The operator
   \begin{align}\label{Dunkl-laplace-drift-gen-hybrid}
 \tilde  L^{\lambda, Dunkl}_k f(x):= &\frac{1}{2}\Delta f(x)  + \sum_{\alpha\in R_+}
k(\alpha)
\frac{\langle\nabla f(x),\alpha\rangle}{\langle\alpha ,x\rangle} + 
\frac{\langle \nabla_x J_k(x,\lambda), \nabla f(x)\rangle }{J_k(x,\lambda)}\\
&-\sum_{\alpha\in R_+}\frac{k(\alpha)\|\alpha\|^2}{2}\cdot 
\frac{f(x)-f(\sigma_\alpha x)}{\langle\alpha , x\rangle^2} \notag
   \end{align}
   for  $f\in C^{(2)}(\mathbb R^N)$ is the generator of a Feller process on $\mathbb R^N$ with the transition densities
   \begin{equation}\label{density-transition-dunkl-drift-hybrid-gen}
 \tilde  K_t^{\lambda,k,Dunkl}(x,A)=\frac{ c_k\cdot e^{-\|\lambda\|^2 t/2}}{t^{\gamma+N/2}} \int_A e^{-(\|x\|^2+\|y\|^2)/(2t)}
 \frac{ E_k(\frac{x}{\sqrt{t}}, \frac{y}{\sqrt{t}}) \cdot J_k(y,\lambda)}{J_k(x,\lambda)}
\cdot w_k(y)\> dy
\end{equation}
for $t>0$, $x\in \mathbb R^N$ and Borel sets $A\subset \mathbb R^N$ with  $ w_k$  as above.
These Feller processes have modifications with RCLL paths and
will be called hybrid Dunkl-Bessel processes of type $R$ with multiplicity $k$ and drift $\lambda$.  
  \end{enumerate}
 \end{def-theorem}

 We now prove part (1) of Theorem \ref{def-dunkl-bessel-p-drift}.
 The proofs of part (2) and of Theorem \ref{def-bessel-p-drift} are completely analog. For this reason we skip the proofs in these cases.

 \begin{proof}[Proof of Theorem   \ref{def-dunkl-bessel-p-drift}(1)]
   In the case $\lambda=0$ without drift this result is shown in \cite{R2} (see Section 2.1 above)
  on the basis of the  reproducing kernel formula
  \begin{equation}\label{reproducing-kernel-dunkl}
 c_k\int_{\mathbb R^N} E_k(y,\lambda) E_k(x,y)e^{-\|y\|^2/2}\cdot w_k(y)\> dy= e^{(\langle x,x\rangle+\langle \lambda, \lambda\rangle)/2} E_k(x,\lambda)
\quad(x,\lambda\in\mathbb C^N)\end{equation}
  from \cite{D1} and some renorming using Proposition \ref{properties-dunkl-kernel}(2). Please notice that our $c_k$ is
   $c_k^{-1}$ in \cite{R2}.
\eqref{reproducing-kernel-dunkl} and Proposition \ref{properties-dunkl-kernel}(2) also yield readily that
$K_t^{\lambda,k,Dunkl}(x,\mathbb R^N)=1$ for all $x\in \mathbb R^N$ and $t>0$, i.e., $(K_t^{\lambda,k,Dunkl})_{t>0}$ is a family of Markov kernels
whose associated transition operators
$$T_t^{\lambda,k,Dunkl}f(x):=\int_{\mathbb R^N} f(y) \> K_t^{\lambda,k,Dunkl}(x,dy)$$
map continuous functions $f\in C_c(\mathbb R^N)$ with compact support into  continuous functions in  $C_0(\mathbb R^N)$,
which are zero at infinity, by the definition of the $K_t^{\lambda,k,Dunkl}$ and by the exponential bounds for $E_k$ in Proposition
\ref{properties-dunkl-kernel}(4). Hence, the  $T_t^{\lambda,k,Dunkl}$ are Feller operators.
Furthermore, \eqref{density-transition-dunkl-drift-gen} implies that
\begin{equation}\label{conjugation-dunkl-1}
T_t^{\lambda,k,Dunkl}f(x) = \frac{e^{-\|\lambda\|^2t/2}}{E_k(x,\lambda)} \cdot T_t^{0,k,Dunkl}(f(.)\cdot E_k(.,\lambda))(x)
\end{equation}
for $t>0$, $x\in\mathbb R^N$, and  $f\in C_b(\mathbb R^N)$. This conjugation identity and the well known semigroup property
$T_t^{0,k,Dunkl}\circ T_s^{0,k,Dunkl}=T_{s+t}^{0,k,Dunkl}$ ($s,t>0$) from \cite{R2} imply that this semigroup property is also available for
$(T_t^{\lambda,k,Dunkl})_{t>0}$.

In order to show that $(T_t^{\lambda,k,Dunkl})_{t\ge 0}$ (with $T_0=Id$) is a Feller semigroup we prove that
\begin{equation}\label{conjugation-dunkl-2}
  \|T_t^{\lambda,k,Dunkl}f -f\|_\infty\to0 \quad\text{for}\quad t\downarrow 0 \quad\text{and all}\quad
  f\in C_c^{(\infty)}(\mathbb R^N),
\end{equation}
which shows that this limit holds for all $f\in C_0(\mathbb R^N)$ by the density of
$C_c^{(\infty)}(\mathbb R^N)$  in $(C_0(\mathbb R^N),\|.\|_\infty$). In order to prove \eqref{conjugation-dunkl-2},
take $f\in  C_c^{(\infty)}(\mathbb R^N)$ and put $g:=E_k(.,\lambda)\cdot f \in C_c^{(\infty)}(\mathbb R^N)$.
Choose constants $0<M<R$ with $supp(f)\subseteq \{x\in\mathbb R^N:\> \|x\|\le M\}$ and put
$m_R:=\inf_{\|x\|\leq R}E_k(x,\lambda)>0$.
For $\|x\|\le R$, we use the fact that $(T_t^{0,k,\mathrm{Dunkl}})_{t\geq0}$ is a Feller semigroup, i.e.,
$$\|T_t^{0,k,\mathrm{Dunkl}}g-g\|_\infty\longrightarrow0
\quad(t\downarrow0).$$
Hence, by \eqref{conjugation-dunkl-1},
$$
T_t^{\lambda,k,\mathrm{Dunkl}}f(x)-f(x)=e^{-\|\lambda\|^2t/2}
\frac{T_t^{0,k,\mathrm{Dunkl}}g(x)-g(x)}{h_\lambda(x)}+(e^{-\|\lambda\|^2t/2}-1)f(x)
$$
and thus
$$\sup_{\|x\|\le R}
\bigl|T_t^{\lambda,k,\mathrm{Dunkl}}f(x)-f(x)\bigr|
\le
\frac{e^{-\|\lambda\|^2t/2}}{m_R}
\|T_t^{0,k,\mathrm{Dunkl}}g-g\|_\infty
+|e^{-\|\lambda\|^2t/2}-1|\|f\|_\infty \to 0
$$
for $t\to0$. Hence,
$$\sup_{\|x\|\le R}
\bigl|T_t^{\lambda,k,\mathrm{Dunkl}}f(x)-f(x)\bigr|\to0.$$
For $\|x\|>R$. we conclude from  Corollary \ref{est-below} that
there are constants $a>0$ and $b\in\mathbb R^N$ such that
$E_k(x,\lambda)\ge a e^{\langle b,x\rangle}$ for $x\in\mathbb R^N$.
Thus,
$$\frac{1}{E_k(x,\lambda)}\le
\frac{ e^{-\langle b,x\rangle}}{a}\le
\frac{e^{\|b\|\|x\|}}{a}.$$
Hence the transition probabilities \eqref{density-transition-dunkl-drift-gen} lead to
\begin{align}
\bigl|T_t^{\lambda,k,\mathrm{Dunkl}}f(x)\bigr|
&\le\frac{c_k e^{-\|\lambda\|^2t/2}}{t^{\gamma+N/2} E_k(x,\lambda)}
\notag\\
&\quad\cdot
\int_{supp(f)}
e^{-(\|x\|^2+\|y\|^2)/(2t)}
E_k(x/\sqrt t,y/\sqrt t)E_k(y,\lambda)|f(y)|w_k(y)\,dy.
\end{align}
Furthermore,  Proposition \ref{properties-dunkl-kernel}(4) yields
$E_k(x/\sqrt t,y/\sqrt t)\le e^{\|x\|\|y\|/t}$, and for $y\in supp(f)$ we have
 $\|y\|\leq M$, and thus for $\|x\|>R$,
$$-\frac{\|x\|^2+\|y\|^2}{2t}
+\frac{\|x\|\|y\|}{t}=
-\frac{(\|x\|-\|y\|)^2}{2t}
\le
-\frac{(\|x\|-M)^2}{2t}.$$
Since
$\int_{supp(f)}
E_k(y,\lambda)||f(y)|w_k(y)\,dy<\infty$,
it follows  that there exists a constant $C>0$
such that for $\|x\|>R$ and $t>0$,
$$
\bigl|T_t^{\lambda,k,\mathrm{Dunkl}}f(x)\bigr|\le
C t^{-\gamma-N/2}\exp\left(
-\frac{(\|x\|-M)^2}{2t}
+\|b\|\|x\|\right).$$
As the  function
$r\longmapsto-\frac{(r-M)^2}{2t}+\|b\|r$
is decreasing on $[R,\infty)$
for sufficiently small $t>0$ we get
$$\sup_{\|x\|\geq R}
\bigl|T_t^{\lambda,k,\mathrm{Dunkl}}f(x)\bigr|
\leq
C t^{-\gamma-N/2}
\exp\left(
-\frac{(R-M)^2}{2t}+\|b\|R
\right)
\to0$$
for $t\to0$ and  $\|x\|>R$.
As $f(x)=0$ for $\|x\|\geq R$, we obtain
$$\sup_{\|x\|\geq R}
\bigl|T_t^{\lambda,k,\mathrm{Dunkl}}f(x)-f(x)\bigr|
\longrightarrow0.$$
This completes the proof of the claim \eqref{conjugation-dunkl-2}.

  In summary, 
$(T_t^{\lambda,k,Dunkl})_{t\ge 0}$ (with $T_0=Id$) is a Feller semigroup.

We now consider  the generators
$$L^{\lambda,Dunkl}_k f:= \lim_{t\downarrow 0}\frac{1}{t}(T_t^{\lambda,k,Dunkl}(f)-f), \quad\quad
L^{0,Dunkl}_k f:= \lim_{t\downarrow 0}\frac{1}{t}(T_t^{0,k,Dunkl}(f)-f)$$
of our transition semigroups with and without drift for
$f\in C_c^{(2)}(\mathbb R^N)$ respectively. Then $f\cdot E_k(.,\lambda)\in C_c^{(2)}(\mathbb R^N)$, and
\begin{align}
  L^{\lambda,Dunkl}_k f(y)&= \lim_{t\downarrow 0}\frac{1}{t}\Biggl(
  \frac{e^{-\|\lambda\|^2t/2}}{E_k(y,\lambda)} \cdot T_t^{0,k,Dunkl}(f(.)\cdot E_k(.,\lambda))(y)-f(y)\Biggr)\notag\\
 &=\lim_{t\downarrow 0}\frac{e^{-\|\lambda\|^2t/2}-1}{t \cdot E_k(y,\lambda)}\cdot T_t^{0,k,Dunkl}(f(.)\cdot E_k(.,\lambda))(y)\notag\\
 &\quad +
  \frac{ 1 }{E_k(y,\lambda)} \lim_{t\to0}
\frac{1}{t}\bigl(   T_t^{0,k,Dunkl}(E_k(.,\lambda) f(.))(y)- E_k(y,\lambda) f(y)\bigr)\notag
\end{align}
and thus
\begin{equation}\label{trafo-generators-dunkl}
 L^{\lambda,Dunkl}_k f(y)= - \frac{\|\lambda\|^2}{2} f(y) +\frac{1}{E_k(y,\lambda)} L^{0,Dunkl}_k (E_k(.,\lambda) f(.))(y).
\end{equation}
$L^{0,Dunkl}_k=\Delta_k/2$,  the product rule in Lemma \ref{prod-rule}, and \eqref{joint-Dunkl-eigenvalues} now yield that
\begin{align}
 L^{\lambda,Dunkl}_k f(y) &=- \frac{\|\lambda\|^2}{2} f(y) + \frac{1}{E_k(y,\lambda)} \Biggl( \frac{f(y)\> \Delta_kE_k(y,\lambda)}{2}
+\frac{\Delta_kf(y)\>E_k(y,\lambda)}{2}\notag\\
&\quad + \langle \nabla f(y), \nabla E_k(y,\lambda)\rangle + 
\sum_{\alpha\in R_+}k(\alpha)\|\alpha\|^2 \frac{(f(y)-f(\sigma_\alpha y))(E_k(y,\lambda)-E_k(\sigma_\alpha y,\lambda))}{2\langle\alpha , y\rangle^2}
\Biggr)\notag\\
&=\frac{\Delta_kf(y)}{2}+  \frac{\langle \nabla f(y), \nabla E_k(y,\lambda)\rangle}{E_k(y,\lambda)}
 +\sum_{\alpha\in R_+}k(\alpha)\frac{\|\alpha\|^2}{2} \Biggl( 1-\frac{ E_k(\sigma_\alpha y,\lambda)}{ E_k(y,\lambda)}\Biggr)
\frac{f(y)-f(\sigma_\alpha y)}{\langle\alpha , y\rangle^2}.
\notag\end{align}
This and \eqref{Dunkl-laplace} now show that $L^{\lambda,Dunkl}_k$ is given as in the theorem.
This completes the proof.
 \end{proof}

 \begin{remark}\label{connection-dunkl-bessel} If one compares the transition probabilities of hybrid Dunkl-Bessel processes
 with drift  in \eqref{density-transition-dunkl-drift-hybrid-gen} with those of 
Bessel processes
 with drift  in \eqref{density-transition-bessel-drift-gen}, it follows by a standard argument on Markov processes (see e.g.
 Exercise (1.17) in Ch. III of \cite{RY}) that the image of a  hybrid Dunkl-Bessel process with drift
 under the canonical projection from
 $\mathbb R^N$ onto the associated Weyl chamber $C$ (which may be regarded as the space of all $W$-orbits of $\mathbb R^N$)
 is just a corresponding Bessel process with drift. In this way, Theorem \ref{def-bessel-p-drift}  is a direct consequence from Theorem
 \ref{def-dunkl-bessel-p-drift}(2).

 We point out that (for $\lambda\ne0$) such a connection is no longer valid between  Dunkl and Bessel processes with drift.
   \end{remark}

 There are several known relations between different Dunkl and Bessel processes without drift; see e.g. \cite{RV1, RV2, CGY, AnV}.
 Most of them hold also for processes with drift. For instance:

\begin{remark}
  For $N\ge2$ consider the root system $R=A_{N-1}$ on $\mathbb R^N$ and  $k\ge0$.
 Let ${\bf 1}:=(1,\ldots,1)\in \mathbb R^N$, and decompose $ \mathbb R^N$ into $\mathbb R\cdot {\bf 1}$
 and its orthogonal complement $$ {\bf 1}^\perp=\{x\in\mathbb R^N:\> \sum_i x_i=0\} \subset \mathbb R^N.$$
Denote the orthogonal projections from  $ \mathbb R^N$ onto  $\mathbb R\cdot {\bf 1}$ and  $ {\bf 1}^\perp$ by
 $\pi_{\bf 1}$ and $\pi_{{\bf 1}^\perp}$ respectively.  Then
\begin{equation}\label{bessel-reduction-A}
 E_k( x,y)= e^{\langle\pi_{\bf 1}(x), \pi_{\bf 1}(y)\rangle}\cdot  E_k( \pi_{{\bf 1}^\perp}(x),
\pi_{{\bf 1}^\perp}(y)) \quad\quad(x,y\in\mathbb R^N).
\end{equation}
Hence, by \eqref{density-transition-dunkl-drift-gen}, for a Dunkl process $(X_t)_{t\ge 0}$ on
$\mathbb R^N$  with drift $\lambda\in \mathbb R^N$, its  orthogonal projection
$$(\pi_{\bf 1}(X_t)=(X_{t,1}+\ldots +X_{t,N})/N)_{t\ge 0}$$ is a one-dimensional Brownian motion on $\mathbb R\cdot {\bf 1}$
with drift $\pi_{\bf 1}(\lambda)=(\lambda_1+\ldots+\lambda_N)/N\cdot {\bf 1}$, and 
$(\pi_{{\bf 1}^\perp}(X_t))_{t\ge 0}$ is a Dunkl process on  $ {\bf 1}^\perp$ (associated with the root system
$\pi_{{\bf 1}^\perp}(R)$ of type $A_{N-1}$ on  $ {\bf 1}^\perp)$ with drift $\pi_{{\bf 1}^\perp}(\lambda)$.
Moreover, these processes are independent whenever this holds for the components
$\pi_{\bf 1}(X_0)$ and  $\pi_{{\bf 1}^\perp}(X_0)$ of $X_0$.
\end{remark}

We next prove that under some (necessary) condition on the initial distributions, the radial parts $(\|X_t\|)_{t\ge0}$ of Dunkl processes 
$(X_t)_{t\ge0}$ with drift on $\mathbb R^N$ are Bessel processes on $[0,\infty[$ with drift,
    and that the same holds for Bessel and hybrid Dunkl-processes with drift. This result generalizes Theorem 4.11 of \cite{RV1}
    for Dunkl processes without drift. This result in \cite{RV1} is based on some connection between the Dunkl kernels $E_k$
    and the one-dimensional spherical Bessel functions
    \begin{equation}\label{def-1dim-besselfkt} j_\alpha(z):=\Gamma(\alpha+1)\sum_{n=0}^\infty \frac{(-1)^n (z/2)^{2n}}{n! \cdot \Gamma(n+\alpha+1)}
  \quad\quad(\alpha\ge -1/2, \> z\in\mathbb C)
    \end{equation}
    which are related to the Bessel functions $J_k(x,y)$ above for $N=1$, $R=\{\pm 1\}$, and $k\ge0$ by
    \begin{equation}\label{connection-bessel-1-dim}
      J_k(x,y)= j_{k-1/2}(ixy); \end{equation}
    see e.g. Eq.~2.1(2) in \cite{RV1} or \cite{D1, R1}.

    We now fix a root system $R$ on $\mathbb R^N$ with some associated closed Weyl chamber $C\subset \mathbb R^N$,
    a multiplicity function $k\ge0$, and a drift vector $\lambda\in\mathbb R^N$. Then, by Eq.~(2.6) in \cite{RV1},
    \begin{equation}\label{rel-dunkl-1-dim-bessel}
      \int_{S^{N-1}} E_k(x,y) w_k(y)\> d\sigma(y)= d_k \cdot j_{\gamma+N/2-1}(i\|x\|)  \quad\quad(x\in\mathbb R^N)
    \end{equation}
    with the surface measure $d\sigma$ on the unit sphere $S^{N-1}\subset  \mathbb R^N$ and the constant
    \begin{equation}\label{connection-ck-dk}
      d_k:= \int_{S^{N-1}}  w_k(y)\> d\sigma(y)=\frac{2}{c_k\cdot \Gamma(\gamma+N/2)}.
      \end{equation}
    In order to describe the necessary condition on the initial distributions, we now say that a probability measure $\mu\in M^1( \mathbb R^N)$
    is Dunkl-$\lambda$-radial if it has the form
     \begin{equation}\label{Dunkl-drift-radial}
\int_{\mathbb R^N} f\> d\mu=\frac{1}{d_k} \int_0^\infty \int_{S^{N-1}} \frac{ f(rz) E_k(\lambda,rz)}{ j_{\gamma+N/2-1}(i\|\lambda\| r)}
 w_k(z)\> d\sigma(z)\> d\nu(r)
     \end{equation}
     for all bounded measurable functions $f$ on $\mathbb R^N$
     and some probability measure $\nu\in M^1([0,\infty[)$. Please notice that the $\mu$  defined in
         \eqref{Dunkl-drift-radial} are in fact probability measures.  With this notation we have:

         \begin{proposition}\label{prp-Dunkl-drift-radial}
           Let  $(X_t)_{t\ge0}$ be a Dunkl process with drift $\lambda\in\mathbb R^N$ on $\mathbb R^N$ such that its initial distribution
           $P_{X_0}$  is Dunkl-$\lambda$-radial. Then for all $t\ge0$, the distributions  $P_{X_t}$  are Dunkl-$\lambda$-radial,
           and  $(\|X_t\|)_{t\ge0}$ is a one-dimensional Bessel process with multiplicity $\gamma+N/2-1/2$ and drift $\|\lambda\|$
in the sense of Theorem \ref{def-bessel-p-drift}.
         \end{proposition}

         Please notice that this result can be applied in particular to the  case $P_{X_0}=\delta_0\in M^1( \mathbb R^N)$.

         \begin{proof}
           Let $\mu:=P_{X_0}\in M^1( \mathbb R^N)$  be related to $\nu\in M^1([0,\infty[)$ as described in  \eqref{Dunkl-drift-radial}.
               Then \eqref{density-transition-dunkl-drift-gen}, \eqref{rel-dunkl-1-dim-bessel}, Fubini, and
               $$ \int_{\mathbb R^N} h(y)\> dy=   \int_0^\infty \int_{S^{N-1}} h(\rho y) \> d\sigma(y) \> \rho^{N-1}\>  d\rho$$
               show that for each
               bounded measurable functions $f$ on $\mathbb R^N$,
               the measure $P_{X_t}$  at time $t>0$ satisfies
               \begin{align} 
                 \int_{\mathbb R^N} f\>& dP_{X_t} = \int_{\mathbb R^N}\int_{\mathbb R^N} f(y) \> K_t^{\lambda,k,Dunkl}(x,dy) \> dP_{X_0}(x)\\
                 &=\frac{ c_k\cdot e^{-\|\lambda\|^2 t/2}}{d_k \> t^{\gamma+N/2}}  \int_0^\infty \int_{S^{N-1}} \Biggl(\int_{\mathbb R^N} 
                 e^{-(r^2+\|y\|^2)/(2t)}f(y) \cdot \frac{E_k(\frac{rz}{\sqrt{t}} \frac{y}{\sqrt{t}}) \cdot E_k(y,\lambda)}{E_k(rz,\lambda)}
                 w_k(y) \> dy\Biggr) \cdot\notag \\
 &\quad\quad\quad\quad\quad\quad\quad\quad\quad\quad\quad\cdot  \frac{E_k(rz,\lambda)}{ j_{\gamma+N/2-1}(i\|\lambda\| r)}
                 w_k(z) d\sigma(z)\> d\nu(r)\notag \\
                  &=\frac{ c_k\cdot e^{-\|\lambda\|^2 t/2}}{d_k \> t^{\gamma+N/2}}  \int_0^\infty \int_{S^{N-1}} \Biggl(\int_{\mathbb R^N} 
                 e^{-(r^2+\|y\|^2)/(2t)}f(y) E_k(\frac{rz}{\sqrt{t}} \frac{y}{\sqrt{t}}) \cdot E_k(y,\lambda)
                 w_k(y) \> dy\Biggr) \cdot\notag \\
 &\quad\quad\quad\quad\quad\quad\quad\quad\quad\quad\quad\cdot  \frac{1}{ j_{\gamma+N/2-1}(i\|\lambda\| r)}
  w_k(z) d\sigma(z)\> d\nu(r)\notag \\
 &= \frac{ c_k\cdot e^{-\|\lambda\|^2 t/2}}{t^{\gamma+N/2}} \int_0^\infty  \int_{\mathbb R^N}  e^{-(r^2+\|y\|^2)/(2t)}f(y)
  \frac{ j_{\gamma+N/2-1}(i r\cdot \|y\|/t)   \cdot E_k(y,\lambda)}{ j_{\gamma+N/2-1}(i\|\lambda\| r)} w_k(y) \> dy
  \> d\nu(r)\notag\\
   &= \frac{  c_k\cdot e^{-\|\lambda\|^2 t/2}}{t^{\gamma+N/2}} \int_0^\infty  \int_0^\infty \int_{S^{N-1}}    e^{-(r^2+\rho^2)/(2t)}f(\rho y)
  \frac{ j_{\gamma+N/2-1}(i r\rho/t)   \cdot E_k(\rho y,\lambda)}{ j_{\gamma+N/2-1}(i\|\lambda\| r)}  w_k(y)\> d\sigma(y)\notag \\
 &\quad\quad\quad\quad\quad\quad\quad\quad\quad\quad\quad\cdot \rho^{2\gamma+N-1}  d\rho
           \> d\nu(r).\notag
               \end{align}
               This shows that $P_{X_t}$ is again  Dunkl-$\lambda$-radial.
               Moreover we conclude from this computation for $\nu=\delta_r$ and from \eqref{rel-dunkl-1-dim-bessel} and \eqref{connection-ck-dk}
               that for all $r\ge0$, $0\le a\le b$,
               $A:=\{x\in\mathbb R^N:\> a\le \|x\|\le b\}$, and $s,t\ge0$, we have  the conditional probabilities
               \begin{align}\label{computation-radial-trans}
                 P(X_{s+t}&\in A| \|X_s\|=r)= \frac{1}{d_k}\int_{S^{N-1}}  \> K_t^{\lambda,k,Dunkl}(rz,A)\frac{   E_k(\lambda,rz)}{ j_{\gamma+N/2-1}(i\|\lambda\| r)}
                 w_k(z)\> d\sigma(z)\\
                 &=\frac{ c_k\cdot e^{-\|\lambda\|^2 t/2}}{t^{\gamma+N/2}} \int_a^b \int_{S^{N-1}}    e^{-(r^2+\rho^2)/(2t)}
  \frac{ j_{\gamma+N/2-1}(i r\rho/t)   \cdot E_k(\rho y,\lambda)}{ j_{\gamma+N/2-1}(i\|\lambda\| r)}  w_k(y)\> d\sigma(y) \rho^{2\gamma+N-1}  d\rho
   \notag\\
                 &=    \frac{d_k c_k\cdot e^{-\|\lambda\|^2 t/2}}{t^{\gamma+N/2}} \int_a^b       e^{-(r^2+\rho^2)/(2t)}    
   \frac{ j_{\gamma+N/2-1}(i r\rho/t) j_{\gamma+N/2-1}(i \rho\cdot\|\lambda\|) }{ j_{\gamma+N/2-1}(i\|\lambda\| r)} \rho^{2\gamma+N-1}  d\rho
    \notag\\
                 &=    \frac{2 e^{-\|\lambda\|^2 t/2}}{\Gamma(\gamma+N/2)\> t^{\gamma+N/2}} \int_a^b       e^{-(r^2+\rho^2)/(2t)}    
  \frac{ j_{\gamma+N/2-1}(i r\rho/t) j_{\gamma+N/2-1}(i \rho\cdot\|\lambda\|) }{ j_{\gamma+N/2-1}(i\|\lambda\| r)} \rho^{2\gamma+N-1}  d\rho.
  \notag\end{align}
  If we use \eqref{connection-bessel-1-dim} and compare the last line of \eqref{computation-radial-trans} with the transition probabilities
  \eqref{density-transition-dunkl-drift-gen} for $N=1$ and multiplicity $\gamma+N/2-1/2$, we see readily that
  $(\|X_t\|)_{t\ge0}$ is a one-dimensional Bessel process with multiplicity $\gamma+N/2-1/2$ and drift $\|\lambda\|$
as claimed.
         \end{proof}

         Proposition \ref{prp-Dunkl-drift-radial} can be transferred to Bessel processes with drift on the Weyl chamber $C$
         as well as to hybrid Dunkl-Bessel-processes with drift on $\mathbb R^N$. For this we have to modify the $\lambda$-radiality of measures in
         \eqref{Dunkl-drift-radial}. We now say that a probability measure $\mu\in M^1(C)$
    is Bessel-$\lambda$-radial if it has the form
     \begin{equation}\label{Bessel-drift-radial}
\int_{C} f\> d\mu=\frac{|W|}{d_k} \int_0^\infty \int_{S^{N-1}\cap C} \frac{ f(rz) J_k(\lambda,rz)}{ j_{\gamma+N/2-1}(i\|\lambda\| r)}
 w_k(z)\> d\sigma(z)\> d\nu(r),
     \end{equation}
     and we say that $\mu\in M^1(\mathbb R^N)$
 is hybrid-Dunkl-Bessel-$\lambda$-radial if it has the form
     \begin{equation}\label{Dunkl-Bessel-drift-radial}
\int_{\mathbb R^N} f\> d\mu=\frac{1}{d_k} \int_0^\infty \int_{S^{N-1}} \frac{ f(rz) J_k(\lambda,rz)}{ j_{\gamma+N/2-1}(i\|\lambda\| r)}
 w_k(z)\> d\sigma(z)\> d\nu(r),
     \end{equation}
     where in both cases $f$ is a bounded measurable function on $C$ or $\mathbb R^N$ and  $\nu\in M^1([0,\infty[)$ is a suitable probability measure.
         Again, measures $\mu$ satisfying \eqref{Bessel-drift-radial} or  \eqref{Dunkl-Bessel-drift-radial} are probability measures.
         The methods of the preceding proposition now lead to the following result where we skip the proof.

\begin{proposition}\label{prp-Dunkl-bessel-drift-radial}
           Let  $(X_t)_{t\ge0}$ be either a Bessel process with drift $\lambda\in C$ on $C$, whose initial distribution
           $P_{X_0}$  is Bessel-$\lambda$-radial,   or a hybrid Dunkl-Bessel process with drift $\lambda\in \mathbb R^N$ on $\mathbb R^N$,
 whose initial distribution
           $P_{X_0}$  is hybrid-Dunkl-Bessel-$\lambda$-radial.
 Then for all $t\ge0$, the distributions  $P_{X_t}$  are Bessel-$\lambda$-radial or  hybrid-Dunkl-Bessel-$\lambda$-radial respectively,
           and  $(\|X_t\|)_{t\ge0}$ is a one-dimensional Bessel process with multiplicity $\gamma+N/2-1/2$ and drift $\|\lambda\|$
in the sense of Theorem \ref{def-bessel-p-drift}.
         \end{proposition}


\section{ A Girsanov-type result}

In this section we present a Girsanov-type result for Dunkl processes, Bessel processes, and hybrid Dunkl-Bessel processes with drift.
This result relates these processes with drift with those without drift. The results and the proofs are analog
to the classical case for Brownian motions
with constant drift. For this we first recapitulate a L\'{e}vy-type characterization for processes without drift; see Proposition 
6.1 of \cite{RV1} in the Dunkl case. The result and its proof there can be easily extended to the following result:

\begin{proposition}\label{char-martingales-processes}
  \begin{enumerate}\itemsep=-1pt
  \item[\rm{(1)}] Let  $(X_t)_{t\ge0}$ be a stochastic process on $\mathbb R^N$ with RCLL paths and start in $0$.
    Then  $(X_t)_{t\ge0}$ is a Dunkl process without drift if and only if
    \begin{equation}\label{exp-mart-dunkl-0}
( E_k(X_t,iz)\cdot e^{t\langle z,z\rangle/2})_{t\ge0}
\end{equation}
    is a $\mathbb C$-valued martingale (w.r.t.~the canonical filtration) for each $z\in\mathbb R^N$. Moreover, the processes 
in \eqref{exp-mart-dunkl-0} are  martingales for all $z\in\mathbb C^N$.
\item[\rm{(2)}] For
  each Dunkl process $(X^\lambda_t)_{t\ge0}$ with drift $\lambda\in\mathbb R^N$ and start in $0$, and for each  $z\in\mathbb C^N$,
 \begin{equation}\label{exp-mart-dunkl-lambda}
\Bigl( \frac{E_k(X_t^\lambda,z)}{E_k(X_t^\lambda,\lambda)}\cdot e^{(-\langle z,z\rangle +\|\lambda\|^2)t/2}\Bigr)_{t\ge0}
\end{equation}
is a martingale.
\item[\rm{(3)}] The statements in (1) and (2) hold also for Bessel processes, when the Dunkl kernels $E_k$
  in \eqref{exp-mart-dunkl-0} and
  \eqref{exp-mart-dunkl-lambda} will be replaced by the Bessel functions $J_k$.
\item[\rm{(4)}] For hybrid Dunkl-Bessel processes $(\tilde X^\lambda_t)_{t\ge0}$ with drift $\lambda\in\mathbb R^N$ and start in $0$, and 
  for each  $z\in\mathbb C^N$,
  \begin{equation}\label{exp-mart-hybrid-dunkl-lambda}
\Bigl( \frac{E_k(\tilde X_t^\lambda,z)}{J_k(\tilde X_t^\lambda,\lambda)}\cdot e^{-(\langle z,z\rangle +\|\lambda\|_2)t/2}\Bigr)_{t\ge0}
\end{equation}
  is a martingale.
  \end{enumerate}
\end{proposition}

We  also need the following classical fact in Girsanov theory; see e.g. \cite{RY}:

\begin{lemma}\label{equiv-measure}
  Let $(\Omega,\cal A,P)$ be a probability space with some filtration, and let $T>0$. Let $(Z_t)_{t\in[0,T]}$ be a $]0,\infty]$-valued martingale
  with $E(Z_0)=1$. Then a $\mathbb C$-valued adapted process $(M_t)_{t\in[0,T]}$ is a martingale w.r.t.~$P$ if and only if $(M_t/Z_t)_{t\in[0,T]}$
  is a martingale w.r.t.~the probability measure $Z_T P$.
  \end{lemma}

With these results we obtain the following Girsanov-type result:

\begin{theorem}\label{girsanov} Let $(X_t^{\lambda,Dunkl})_{t\ge0}$ be a Dunkl process on $\mathbb R^N$ with drift $\lambda\in\mathbb R^N$
  and start in $0\in\mathbb R^N$, which is defined on some probability space $(\Omega,\cal A, P)$.
  Then for each $T>0$,  $(X_t^{\lambda,Dunkl})_{t\in[0,T]}$ is a Dunkl process without drift with start in $0\in\mathbb R^N$
  w.r.t.~the  probability measure $\frac{e^{\|\lambda\|^2 T/2}}{E_k(X_T,\lambda)}\cdot P$ on  $(\Omega,\cal A)$.
\end{theorem}

\begin{proof}
  Let $z\in\mathbb R^N$. Then, by Proposition \ref{char-martingales-processes}(2),
  $$\Bigl( \frac{E_k(X_t,iz)e^{\|z\|^2 t/2} }{E_k(X_t,\lambda)e^{-\|\lambda\|^2 t/2}}\Bigr)_{t\ge0}$$
  is a martingale w.r.t.~$P$. In particular, for $z=0$,
  $$Z_t:= \Bigl(\frac{1}{E_k(X_t,\lambda)e^{-\|\lambda\|^2 t/2}}\Bigr)_{t\ge0}$$ is a positive martingale  w.r.t.~$P$.
  Hence, by Lemma \ref{equiv-measure},
  $( E_k(X_t,iz)e^{\|z\|^2 t/2} )_{t\in[0,T]}$ is a martingale  w.r.t.~$Z_T P$ for each  $z\in\mathbb R^N$.
    The claim now follows from  Proposition \ref{char-martingales-processes}(1).
\end{proof}

\begin{remark} The assertion of Theorem \ref{girsanov} holds also for Bessel processes and hybrid Dunkl-Bessel processes
where here in both cases the Dunkl kernel $E_k$ has to be replaced by the Bessel function $J_k$ in the density by
parts (3) and (4) of  Proposition \ref{char-martingales-processes}.
\end{remark}

\begin{remark}
The generators \eqref{Bessel-laplace-drift-gen}, \eqref{Dunkl-laplace-drift-gen}, and \eqref{Dunkl-laplace-drift-gen-hybrid}
of Bessel, Dunkl, and hybrid Dunkl-Bessel processes with drift respectively
lead to stochastic differential equations (SDEs) for these processes.
For these processes without drift, these SDEs are studied in particular in \cite{CGY, GY1, GY2, GY3}. For general existence and uniqueness of
associated initial value problems we refer to \cite{GM}. These existence and uniqueness results hold also for the processes
with drift on the compact time intervals  $[0,T]$ ($T>0$), as
the Girsanov-type theorem \ref{girsanov} transfers solutions without drift into
solutions with drift and vice versa.
  \end{remark}

\section{Moments, strong laws of large numbers, and central limit theorems for Dunkl processes with drift}

For the multiplicity $k=0$, Dunkl processes on $\mathbb R^N$ without drift are Brownian motions, and, by \cite{R2, RV1, RV2},
Dunkl processes  without drift with $k\ge0$ have much in common with
Brownian motions on $\mathbb R^N$. In particular, they are martingales (in all coordinates), there exist associated heat polynomials, and
a stochastic analysis  as well as a Wiener chaos can be established; see  \cite{CGY, GY1, GY2, GY3}.
We show that Dunkl processes on $\mathbb R^N$ with drift $\lambda\in\mathbb R^N$ also have much in common with
Brownian motions on $\mathbb R^N$ with drift $\lambda$.

For instance, we have the following result for Dunkl processes $(X_t^\lambda)_{t\ge0}$
with drift $\lambda\in\mathbb R^N$.

\begin{theorem}\label{expectation-dunkl-drift} Let $(X_t^\lambda)_{t\ge0}$ be a  Dunkl process
  with drift $\lambda\in\mathbb R^N$ and $k\ge0$ such that the starting distribution  $P_{X_0^\lambda}\in M^1(\mathbb R^N)$ has first moments.
  Then $(X_t^\lambda-t\lambda)_{t\ge0}$ is a martingale w.r.t.~the canonical filtration
of  $(X_t^\lambda)_{t\ge0}$.
In particular, for $P_{X_0^\lambda}=\delta_0$, $E(X_t^\lambda)=t\lambda$.
\end{theorem}

The proof of this theorem is based on the following observation:

\begin{lemma}\label{expectation-dunkl-drift-xi}  Let $\lambda\in\mathbb R^N$, $i=1,\ldots,N$, and $p_i(x):=x_i$ ($x\in\mathbb R^N$).
  Then $L^{\lambda, Dunkl}_k (p_i)(x)=\lambda_i$.
\end{lemma}

\begin{proof} We conclude from \eqref{Dunkl-laplace-drift-gen} that
\begin{align}
L_k^{\lambda, Dunkl}(p_i)(x)=& \sum_{\alpha\in R_+}k(\alpha)
\frac{\alpha_i}{\langle\alpha ,x\rangle} + 
\frac{\partial_{x_i}E_k(x,\lambda)}{E_k(x,\lambda)}\notag\\
&\quad -\sum_{\alpha\in R_+}  \frac{k(\alpha)\|\alpha\|^2}{2} \frac{E_k(\sigma_\alpha x,\lambda)}{E_k(x,\lambda)}\cdot
\frac{p_i(x)-p_i(\sigma_\alpha x)}{\langle\alpha ,x\rangle^2}.\notag
\end{align}
Using $(\sigma_\alpha x)_i=x_i-2\frac{\langle\alpha ,x\rangle}{\|\alpha\|^2}\alpha_i$ and \eqref{joint-Dunkl-eigenvalues}, we thus obtain
\begin{align}
  L_k^{\lambda, Dunkl}(p_i)(x)=&\frac{1}{E_k(x,\lambda)}\Biggl( \partial_{x_i}E_k(x,\lambda)+
  \sum_{\alpha\in R_+}k(\alpha)\frac{\alpha_i}{\langle\alpha ,x\rangle}E_k(x,\lambda)
  -\sum_{\alpha\in R_+}k(\alpha)\>  E_k(\sigma_\alpha x,\lambda) \frac{\alpha_i}{\langle\alpha ,x\rangle}\Biggr)\notag\\
  =&\frac{1}{E_k(x,\lambda)}\Biggl( \partial_{x_i}E_k(x,\lambda)+\sum_{\alpha\in R_+}\frac{k(\alpha)}{\langle\alpha ,x\rangle}
  (E_k(x,\lambda)-E_k(\sigma_\alpha x,\lambda) )\alpha_i\Biggr)\notag\\
  =& \frac{T_{e_i}E_k(x,\lambda)}{E_k(x,\lambda)}\quad=\quad \lambda_i\notag
\end{align}
(where $T_{e_i}$ stands for the Dunkl operator  $T_{e_i}(k)$)
which proves the  statement.
\end{proof}

\begin{proof}[Proof of Theorem \ref{expectation-dunkl-drift}]
  $L_k^{\lambda, Dunkl}c=0$ for constants $c$,  Lemma \ref{expectation-dunkl-drift-xi}, and Theorem \ref{def-dunkl-bessel-p-drift}
  show that
  $$ \int_{\mathbb R^N}y_i  \> K_t^{\lambda,k,Dunkl}(x,dy)  =     e^{t\cdot L_k^{\lambda, Dunkl}}(p_i)(x)=x_i+t\lambda_i \quad\quad (t\ge 0, \> x\in\mathbb R^N) $$
  and thus $E(X_{s+t,i}-(s+t)\lambda_i| \cal F_s)=X_{s,i}-s\lambda_i$ for $s,t\ge0$ and the canonical filtration  $(\cal F_t)_{t\ge0}$
  as claimed.
\end{proof}

\begin{remark} Having corresponding results for higher moments of Dunkl processes without drift in \cite{RV1, RV2} in mind,
  one might expect that Theorem \ref{expectation-dunkl-drift} and Lemma \ref{expectation-dunkl-drift-xi}
  can be  extended to higher moments. However, for $i\ne j$ and $p_{i,j}(x):=x_ix_j$,  the methods in the proof of Lemma \ref{expectation-dunkl-drift-xi} yield
\begin{align}
L_k^{\lambda, Dunkl}(p_{i,j})(x)=& \sum_{\alpha\in R_+}k(\alpha)
\frac{\alpha_i x_j+ \alpha_j x_i}{\langle\alpha ,x\rangle} + 
\frac{(x_j\partial_{x_i}+x_i\partial_{x_j})E_k(x,\lambda)}{E_k(x,\lambda)}\notag\\
&\quad -\sum_{\alpha\in R_+}  \frac{k(\alpha)\|\alpha\|^2}{2} \frac{E_k(\sigma_\alpha x,\lambda)}{E_k(x,\lambda)}\cdot
\frac{p_{i,j}(x)-p_{i,j}(\sigma_\alpha x)}{\langle\alpha ,x\rangle^2}.\notag
\end{align}
As
$$p_{i,j}(x)-p_{i,j}(\sigma_\alpha x)=2 \frac{\langle\alpha ,x\rangle}{\|\alpha\|^2}(x_j\alpha_i+x_i\alpha_j)-
4 \frac{\langle\alpha ,x\rangle^2}{\|\alpha\|^4}\alpha_i\alpha_j,$$
we obtain 
\begin{equation}\label{highre-moments-trial}
  L_k^{\lambda,  Dunkl}(p_{i,j})(x)=
  \lambda_i x_j+\lambda_j x_i  + 2\sum_{\alpha\in R_+} \frac{k(\alpha)}{\|\alpha\|^2}\alpha_i\alpha_j\frac{E_k(\sigma_\alpha x,\lambda)}{E_k(x,\lambda)}.
\end{equation}
Moreover, for $i=j$ and $p_{i,i}(x):=x_i^2$ and for higher moments we obtain a similar result with the same sum
in the end of the right hand side such that the nice approach above for the first moments does not seem to work for higher moments of Dunkl
processes with drift $\lambda\ne0$.

  The Dunkl operators $T_\xi(k)$  map polynomials (in $N$ variables) of degree $n\ge 1$ into polynomials of degree $n-1$.
  This might  suggest that conjugations
  $$T_\xi^\lambda(k)(f)(x):= \frac{1}{E_k(x,\lambda)} T_\xi(k)(E(.,\lambda)\cdot f)(x)$$
  with drift vectors $\lambda$  also map polynomials into polynomials. This is  not correct in general. For instance,
  for $N=1$, $\lambda=1$, $\xi=1$, $R_+=\{1\}$, and $f(x)=x$, a short computation yields
  $T_\xi^\lambda(k)(f)(x)=1+x +2k\frac{E_k(-x,1)}{E_k(x,1)}$, which usually fails to be a polynomial; see e.g. Example \ref{example-11}.
  \end{remark}

In summary, usual higher moments do not seem to be adapted to Dunkl processes with drift $\lambda\ne0$.
It can be checked that this also holds for Bessel processes and hybrid Dunkl-Bessel processes with drift $\lambda\ne 0$
(even for moment of first order). In this way,
 Theorem \ref{expectation-dunkl-drift} and Lemma  \ref{expectation-dunkl-drift-xi} are remarkable singular results.
In order to obtain strong laws of large numbers for our processes with drift,
we now use the approach of the theory of modified moments for random walks on commutative hypergroups
in \cite{Z2}, Ch.~7 of \cite{BH}, and references cited there where we restrict our attention to first and second moments for practical reasons.

We now fix $R,R_+, W, C, k\ge0, E_k$, and the drift vector
$\lambda\in\mathbb R^N$, and study this approach for the generators $L_k^{\lambda, k, Dunkl}$, i.e. Dunkl processes with drift, first.

\begin{definition}\label{def-moments-dunkl} Define the
  multivariate (modified)
  moment functions of order 1 and 2 by
  $$m_1^{Dunkl}(x):=(m_{1;1}^{Dunkl}(x),\ldots,m_{1;N}^{Dunkl}(x)):= \frac{1}{E_k(x,\lambda)}\nabla_\lambda E_k(x,\lambda)$$
  and
  $$m_2^{Dunkl}(x):=(m_{2;j,l}^{Dunkl}(x))_{j,l=1,\ldots,N}:= \frac{1}{E_k(x,\lambda)}H_\lambda E_k(x,\lambda)$$
  for $x\in\mathbb R^N$ with the gradient $\nabla_\lambda$ and the Hessian $H_\lambda$ w.r.t.~the variable
  $\lambda\in\mathbb R^N$ respectively.
\end{definition}

For abbreviation, we now suppress the superscript $Dunkl$ for these moment functions.

The integral representation \ref{general-positive-integralrep-dunkl}  and the methods in Section 7.2 of \cite{BH} lead to:

\begin{lemma}\label{bounds-moments-dunkl}
  For all $l,j=1,\ldots,N$ and $x\in\mathbb R^N$:
\begin{enumerate}\itemsep=-1pt
  \item[\rm{(1)}] $m_{1;l}(x)^2\le m_{2;l,l}(x).$
  \item[\rm{(2)}] $m_{2;l,j}(x)^2\le m_{2;l,l}(x)\cdot m_{2;j,j}(x).$
\item[\rm{(3)}] $ \sum_{l=1}^N  m_{2;l,l}(x) \le \|x\|^2.$
\item[\rm{(4)}] $\|m_1(x)\| \le \|x\|.$
  \item[\rm{(5)}]  $m_{1;l}(0)=m_{2;l,j}(0)=0$.
 \end{enumerate}
\end{lemma}

\begin{proof} (1) follows from Theorem \ref{general-positive-integralrep-dunkl} and Jensen's inequality,
  while for (2) the Cauchy-Schwarz inequality has to be used. (3)  follows from Theorem  \ref{general-positive-integralrep-dunkl}(1).
  (4) follows from (1) and (3), and (5)  from $E_k(\lambda,0)=1$
 for  $\lambda\in\mathbb R^N$.
\end{proof}

We next show that $m_1, m_2$  are adapted  the generators $L_k^{\lambda, Dunkl}$ as follows:

\begin{proposition}\label{mart-m1} For all $x,\lambda\in\mathbb R^N$,
  $$L_k^{\lambda,  Dunkl}m_1(x)=\lambda.$$
Moreover, if $(X_t^\lambda)_{t\ge0}$ is a  Dunkl process
  with drift $\lambda$  such that the starting distribution  $P_{X_0^\lambda}\in M^1(\mathbb R^N)$ has first moments,
  then $(m_1(X_t^\lambda)-t\lambda)_{t\ge0}$ is a martingale.
In particular, for $P_{X_0^\lambda}=\delta_0$, $E(m_1(X_t^\lambda))=t\lambda$.
\end{proposition}

\begin{proof} Eq.~\eqref{trafo-generators-dunkl} yields that
\begin{align}\label{generator-m1} 
  L_k^{\lambda,Dunkl}m_1(x)=& - \frac{\|\lambda\|^2}{2} m_1(x)
  +\frac{1}{E_k(x,\lambda)} L^{0,k,Dunkl}(\nabla_\lambda E_k(x,\lambda)) \notag\\
  =& - \frac{\|\lambda\|^2}{2} m_1(x)
  +\frac{1}{E_k(x,\lambda)} \nabla_\lambda( L^{0,k,Dunkl}E_k(x,\lambda)) \notag\\
   =& - \frac{\|\lambda\|^2}{2} m_1(x)
   +\frac{1}{2 E_k(x,\lambda)} \nabla_\lambda( \|\lambda\|^2E_k(x,\lambda)) \notag\\
   =& - \frac{\|\lambda\|^2}{2} m_1(x)+ \frac{\|\lambda\|^2}{2}\frac{\nabla_\lambda E_k(x,\lambda)}{E_k(x,\lambda)}
   +\lambda \quad =\quad \lambda
\end{align}
as claimed. Please notice that we used the known fact 
for the second $=$ that $\nabla_\lambda$ and   $L_k^{0,Dunkl}$ (applied w.r.t.~$x$)
commute. As in the proof  of Theorem \ref{expectation-dunkl-drift} we thus obtain
$$\int_{\mathbb R^N} m_1(y)\> K_t^{\lambda, k,Dunkl}(x,dy)= m_1(x)+t\lambda  \quad\quad(x\in\mathbb R^N)$$
and hence the remaining statements   where $E(m_1(X_t^\lambda))$
exists for all $t\ge0$ due to Lemma \ref{bounds-moments-dunkl}(4).
\end{proof}

\begin{remark} Proposition \ref{mart-m1}, Theorem \ref{expectation-dunkl-drift}, and Lemma \ref{expectation-dunkl-drift-xi}  imply that
  all components of  the function $x\mapsto m_1(x)-x$ are $L_k^{\lambda,  Dunkl}$-harmonic, and that
  the processes $( m_1(X_t^\lambda)-X_t^\lambda)_{t\ge0}$ are martingales. For $k=0$ or $\lambda=0$, this is trivial (as here $m_1(x)=x$),
  but in all other cases
  this observation is  remarkable. We point out that such a result  never holds for the Sturm-Liouville hypergroups on $[0,\infty[$
      studied in Ch.~7 of \cite{BH}. Also in the case of Bessel and hybrid Dunkl-Bessel processes, such a result is not available.
\end{remark}

We now consider the second modified moments.

\begin{proposition}\label{mart-m2}  For all $x,\lambda\in\mathbb R^N$ and $j,l=1,\ldots,N$, 
  $$  L_k^{\lambda,Dunkl}m_{2; j,l}(x) =\delta_{j,l} + \lambda_j m_{1,l}(x) +\lambda_l m_{1,j}(x),$$
Moreover, if $(X_t^\lambda)_{t\ge0}$ is a  Dunkl process
  with drift $\lambda$  such that the starting distribution  $P_{X_0}^\lambda\in M^1(\mathbb R^N)$ has second moments,
  then for all $j,l=1,\ldots,N$,
  \begin{align}\label{mart-m2-cent}
  ( m_{2; j,l}(X_t^\lambda)- E( m_{2; j,l}(X_t^\lambda)) &- m_{1;l}(X_t^\lambda)E( m_{1;j}(X_t^\lambda))- m_{1;j}(X_t^\lambda)E( m_{1;l}(X_t^\lambda))
  \notag \\ &+2E( m_{1;l}(X_t^\lambda))E( m_{1;j}(X_t^\lambda)) )_{t\ge0}
  \end{align}
  is a martingale.
  In particular, for $P_{X_0}^\lambda=\delta_0$ and $l,j=1,\ldots,N$,
  \begin{equation}\label{2nd-mod-moments}
    E(m_{2; l,j}(X_t^\lambda))= t\delta_{j,l}+t^2 \lambda_j\lambda_l.\end{equation}
\end{proposition}

\begin{proof} As in the preceding proof, Eq.~\eqref{trafo-generators-dunkl} yields for $j,l=1,\ldots,N$ that
\begin{align}\label{generator-m2} 
  L_k^{\lambda,Dunkl}m_{2; j,l}(x)=& - \frac{\|\lambda\|^2}{2} m_{2; j,l}(x)
  +\frac{1}{E_k(x,\lambda)} L^{0,k,Dunkl}(\partial_{\lambda_j}\partial_{\lambda_l} E_k(x,\lambda)) \notag\\
  =& - \frac{\|\lambda\|^2}{2}  m_{2; j,l}(x)
  +\frac{1}{E_k(x,\lambda)}\partial_{\lambda_j}\partial_{\lambda_l} ( L^{0,k,Dunkl}E_k(x,\lambda)) \notag\\
   =& - \frac{\|\lambda\|^2}{2} m_{2; j,l}(x)
   +\frac{1}{2 E_k(x,\lambda)} \partial_{\lambda_j}\partial_{\lambda_l}( \|\lambda\|^2E_k(x,\lambda)) \notag\\
   =& \delta_{j,l} + \lambda_j m_{1;l}(x) +\lambda_l m_{1;j}(x)
\end{align}
which proves the first claim.
Therefore, as $L_k^{\lambda,Dunkl}c=0$ for constants $c$, and by \eqref{generator-m1}, 
\begin{align}\label{generator-m2-2}
  \int_{\mathbb R^N} m_{2;l,j}(y) \> K_t^{\lambda,k,Dunkl}(x,dy)  =&     e^{t\cdot L_k^{\lambda,  Dunkl}}(m_{2;j,l})(x) \\
   =&   m_{2;l,j}(x)  + t   L_k^{\lambda,Dunkl}m_{2; j,l}(x)+ \frac{t^2}{2} (   L_k^{\lambda,Dunkl})^2 m_{2; j,l}(x)+\ldots \notag\\
  =&  m_{2;j,l}(x) + t( \delta_{j,l} +\lambda_l  m_{1;j}(x)+\lambda_j  m_{1;l}(x)) + t^2\lambda_l\lambda_j
\notag\end{align}
and thus, for $s,t\ge0$,
$$E(m_{2;l,j}(X_{s+t}^\lambda))=E(m_{2;l,j}(X_{s}^\lambda))+t( \delta_{j,l} +\lambda_l  E(m_{1;j}(X_{s}^\lambda)) +\lambda_j  E(m_{1;l}(X_{s}^\lambda)) ) + t^2\lambda_l\lambda_j$$
where all expectations exist by Lemma \ref{bounds-moments-dunkl}. This and Lemma \ref{bounds-moments-dunkl}(5)
lead to the last identity in the proposition. Finally, we see that for $s,t\ge0$ and the canonical filtration  $(\cal F_t)_{t\ge0}$,
$$ E(m_{2;l,j}(X_{s+t}^\lambda)| \cal F_s)=m_{2;l,j}(X_{s}^\lambda)+ t(\delta_{j,l} +\lambda_l  m_{1;j}(X_{s}^\lambda)+\lambda_j  m_{1;l}(X_{s}^\lambda) )+ t^2\lambda_l\lambda_j.$$
As  $E(m_{1;j}(X_{s+t}^\lambda)| \cal F_s)=m_{1;j}(X_{s}^\lambda) +t\lambda_j$ by Proposition \ref{mart-m1}, it follows by a short computation that
the process in \eqref{mart-m2-cent} is in fact a martingale as claimed.
\end{proof}

The preceding results on the moment functions can be used to derive a strong law of large numbers
similar to Theorem 7.3.21 in \cite{BH}. However, as we deal with processes in continuous time here, the details are slightly more involved here. We need the following two  martingale results  in continuous time.

\begin{theorem}\label{mart-convergence-SLLN}
  \begin{enumerate}\itemsep=-1pt
  \item[\rm{(1)}] Let $(X_t)_{t\ge0}$ be a submartingale with a.s.~RCLL paths and with
    $$\sup_{t\ge0} E(max(X_t,0))< \infty.$$ Then  
    $\lim_{t\to\infty} X_t$ exists a.s.~as an integrable random variable.
  \item[\rm{(2)}] {\bf Continuous martingale Kronecker lemma:}
    Let $(X_t)_{t\ge0}$ be a martingale with a.s.~RCLL paths and $a:[0,\infty[\to]0,\infty[$ a continuous,
        increasing function with $\lim_{t\to\infty} a(t)=\infty$.
        If  $\lim_{t\to\infty} \int_0^t \frac{1}{a(s)} \> dX_s$ exists a.s (in $\mathbb R$), then  $\lim_{t\to\infty}  X_t/a(t)=0$ a.s..
 \end{enumerate}
\end{theorem}

\begin{proof} (1) is well-known; see e.g.~\cite{LS}.
  In fact, Doob's upcrossing inequality implies that the statement for $t\in\mathbb Q$,
  and the general case follows from the RCLL property. The statement in (2) appears in the literature under
  slightly different conditions; see e.g.~Section 4.2 of \cite{LS} or \cite{EM}.
  An analysis of the proofs there shows that the slightly different and usually much more general conditions there
  are used only for the ``product formula''
  $$ \frac{X_t}{a(t)} -  \frac{X_0}{a(0)}  - \int_0^t X_s \> d  \frac{1}{a(s)} - \int_0^t \frac{1}{a(s)} \> dX_s=[ \frac{1}{a}, X]_t=0$$
with the mutual quadratic variation $[ ., .]$ which holds under our assumptions. We therefore skip the details of the proof of (2) here.
\end{proof}

We now apply Theorem \ref{mart-convergence-SLLN} to the following strong law:

\begin{proposition}\label{SLLN-m1}
   Consider a Dunkl process  $(X_t^\lambda)_{t\ge0}$ with RCLL paths,
    drift $\lambda$,  and start in 0. Then for all $\epsilon>1/2$,
   $$\lim_{t\to\infty} \frac{m_1(X_t^\lambda)-t\lambda}{t^\epsilon}=0 \quad{a.s.}.$$
\end{proposition}

\begin{proof} For $j=1,\ldots,N$ consider the process
  $$(V_{t,j}:= m_{2;j,j}(X_t^\lambda)-
  2 m_{1;j}(X_t^\lambda)E( m_{1;j}(X_t^\lambda))+E( m_{1;j}(X_t^\lambda))^2)_{t\ge0}.$$
We know from Propositions \ref{mart-m1} and \ref{mart-m2} that
$$ V_{t,j}=m_{2;j,j}(X_t^\lambda) -2t\lambda_j  m_{1;j}(X_t^\lambda) +t^2\lambda_j^2 \quad\quad (t\ge0),$$
and that $(V_{t,j}-t)_{t\ge0}$ is a martingale. Furthermore,  by Lemma \ref{bounds-moments-dunkl}(1),
\begin{equation}\label{est-pos-v}
0\le(m_{1;j}(X_t^\lambda)-t\lambda_j)^2\le V_{t,j} \quad\quad (t\ge0),
  \end{equation}
i.e., $(V_{t,j})_{t\ge0}$ is a nonnegative RCLL submartingale. Hence,
$$\Bigl(\tilde  V_{t,j}:= \int_0^t \frac{1}{1+s^{2\epsilon}}\> d V_{s,j}=  \int_0^t \frac{1}{1+s^{2\epsilon}}\> d( V_{s,j}-s)+
\int_0^t \frac{1}{1+s^{2\epsilon}}\> ds\Bigr)_{t\ge0}$$
is a nonnegative RCLL submartingale with
$$\lim_{t\to\infty} E(\tilde  V_{t,j}) = \int_0^\infty  \frac{1}{1+s^{2\epsilon}}\>ds<\infty.$$
Hence, $\tilde  V_{t,j}$ converges for $t\to\infty$ to some integrable random variable a.s. by Theorem
\ref{mart-convergence-SLLN}(1). This means that the martingale 
$$\Bigl(  \int_0^t\frac{1}{1+s^{2\epsilon}}\> d( V_{s,j}-s)\Bigr)_{t\ge0}$$
also converges a.s.. Hence, by Theorem
\ref{mart-convergence-SLLN}(2), 
$$\lim_{t\to\infty} \frac{1}{1+t^{2\epsilon}} (V_{t,j}-t-V_{0,j})=0\quad{a.s.}$$
and thus $\lim_{t\to\infty} V_{t,j}/t^{2\epsilon}=0$ a.s.. The inequality
\eqref{est-pos-v} now completes the proof.
\end{proof}

In the next step we try to remove the first moment function $m_1$ in the statement of Proposition \ref{SLLN-m1}
(like in the case of Sturm-Liouville hypergroups in Ch.~7 of \cite{BH}). For this we consider the projection $p_C:\mathbb R^N\to C$
onto some Weyl chamber $C$  which assigns to each $x\in\mathbb R^N$ its Weyl group orbit, where the space of Weyl group orbits is identified with 
$C$.

\begin{lemma}\label{limit-m1-Dunkl} Let $k>0$.
  Assume that $\lambda\in\mathbb R^N$ is not contained in the boundary of any Weyl chamber, i.e., 
  $\lambda$ is contained in some unique Weyl chamber $C\subset\mathbb R^N$. Then, for each $x\in\mathbb R^N$,
  $$\lim_{t\to\infty} \frac{m_1(tx)-p_C(tx)}{t}=0.$$
\end{lemma}

\begin{proof} Fix $j=1,\ldots,N$. By the integral representation \ref{general-positive-integralrep-dunkl}
  $$m_{1;j}(tx)= \frac{\partial_{\lambda_j} E_k(tx,\lambda)}{E_k(tx,\lambda)}=
  t\frac{ \int_{\mathbb R^N} y_j e^{t\langle y,\lambda\rangle}\> d\mu_x(y)}{\int_{\mathbb R^N} e^{t\langle y,\lambda\rangle}\> d\mu_x(y)} $$
  with  $\{gx: \> g\in W\}\subset  supp\>\mu_x\subset co(\{gx: \> g\in W\})$ where  $co$ stands for the convex hull.
  As $\langle y,\lambda\rangle$ for $y\in\{gx: \> g\in W\}$ is maximal for $y=p_C(x)$ (this follows e.g. from the lemma in Section 1,12 of \cite{H}),
  the linear function $y\mapsto \langle y,\lambda\rangle$ takes its maximum on   $supp\>\mu_x$ precisely in $p_C(x)$.
  This implies that the probability measures
  $$\frac{1}{\int_{\mathbb R^N} e^{t\langle y,\lambda\rangle}\> d\mu_x(y)}  e^{t\langle y,\lambda\rangle}\> d\mu_x(y)$$
  tend weakly to $\delta_{p_C(x)}$ for $t\to\infty$ where these measures have a fixed compact support. This yields the lemma.
\end{proof}

For $N=1$ (and $R=\{\pm 1\}$), Lemma \ref{limit-m1-Dunkl} trivially implies that for $k>0$,
\begin{equation}\label{complete-limit-m1-1dim}
  \lim_{|x|\to\infty} \frac{m_1(x)}{|x|}= sgn(\lambda).\end{equation}
  This
and Proposition \ref{SLLN-m1} yield:

\begin{corollary}\label{complete-slln-1dim-Dunkl}
  For a $1$-dimensional Dunkl process  $(X_t^\lambda)_{t\ge0}$ with RCLL paths, $k>0$,
    drift $\lambda\ne0$,  and start in 0,
   $$\lim_{t\to\infty} \frac{|X_t^\lambda|}{t} = |\lambda| \quad{a.s.}.$$
\end{corollary}

We next study the jump rates of the generator $L_k^{\lambda,Dunkl}$ in \eqref{Dunkl-laplace-drift-gen} for $N=1$. This will lead to:

\begin{theorem}\label{SLLN-1dim-dunkl}
  For a $1$-dimensional Dunkl process  $(X_t^\lambda)_{t\ge0}$ with RCLL paths, $k>0$,
    drift $\lambda\ne0$,  and start in 0,
   $$\lim_{t\to\infty} \frac{X_t^\lambda}{t} = \lambda \quad{a.s.}.$$
\end{theorem}

\begin{proof} It is sufficient to prove the result for $\lambda>0$. It is known (see e.g. \cite{R1}) that for $N=1$,
  $$E_k(x,\lambda)= e^{x \lambda}\>_1F_1(k,2k+1; -2x \lambda).$$
  Moreover, using the Kummer transformation for $_1F_1$ (see e.g.~13.2.39 of \cite{NIST}) as well as the limit relation for 
  $_1F_1(a,b;z)$ for $z\to\infty$ in 13.2.23 of \cite{NIST}, we have
  \begin{align}\label{limit-1f1}
    E_k(x,\lambda)\sim& \frac{\Gamma(2k+1)}{\Gamma(k+1)2^k} e^{x \lambda}\cdot (x\lambda)^{-k} \quad\quad\text{for}\quad\quad x\to\infty
 \quad\quad\text{and}   \notag\\
 E_k(x,\lambda)\sim& \frac{\Gamma(2k+1)}{\Gamma(k)2^{k+1}} e^{-x \lambda}\cdot (-x\lambda)^{-(k+1)}
 \quad\quad\text{for}\quad\quad x\to-\infty.\end{align}
  Therefore, the jump rates $R(x)$ of $L_k^{\lambda,Dunkl}$  in \eqref{Dunkl-laplace-drift-gen} from $x\ne0$ to $-x$ satisfy
  \begin{equation} R(x)\sim \frac{k^2}{4\lambda} \cdot \frac{1}{x^3} \quad\text{for}\quad x\to\infty
 \quad\quad\text{and} \quad\quad R(-x)\sim \frac{\lambda}{x}  \quad\text{for}\quad x\to-\infty.
    \end{equation}
As $\int_1^\infty x^{-3}\> dx <\infty $ and $\int_1^\infty  x^{-1}\> dx =\infty $, we conclude from Corollary
\ref{complete-slln-1dim-Dunkl} that, almost surely, the number of jumps from $]0,\infty[$ to $]-\infty,0[$
        for $t>0$ sufficiently large is finite, while the number of jumps  from $]-\infty,0[$ to $]0,\infty[$ is infinite on
         $[t_0,\infty[$ for $t_0>0$    under the condition that there are infinitely many jumps  from $]0,\infty[$ to $]-\infty,0[$.
        This yields the claim.
      \end{proof}

\begin{remark} A Dunkl process
  $(X_t^\lambda)_{t\ge0}$ with $k=0$ and drift $\lambda$ is just a usual Brownian motion with drift $\lambda$,
  and the moment function $m_1$ above is  $m_1(x)=x$. This means that  here,
  \begin{equation}\label{conjecture-SLLN}\lim_{t\to\infty} \frac{X_t^\lambda}{t} =\lambda \quad{a.s.}.\end{equation}
  We conjecture  that this statement also holds for all $k>0$, $N\ge2$,  all root systems, and all drift vectors $\lambda$
  (possibly under the condition of not being on the boundary of some Weyl chamber).
As a first step for the proof, we would have to show that Lemma \ref{limit-m1-Dunkl} 
  can be 
  improved to the  statement
  \begin{equation}\label{conjecture-m1} \lim_{\|x\|\to\infty} \frac{m_1(x)-p_C(x)}{\|x\|}= 0\end{equation}
    like for $N=1$.  Unfortunately, this is not clear, as for $N\ge2$ not
    much information about the representing measures $\mu_x$ are known
    (except for the  case where the Weyl group is $\mathbb Z^N$, where this follows from the one-dimensional case).
    In this context we point out that
    there exist interesting particular limits and estimates for $E_k(x,\lambda)$ and the Dunkl heat kernels
in \eqref{density-transition-dunkl-gen}
for $\|x\|\to\infty$; see e.g. \cite{RdJ, AG, GS, AT, DH} and references there.
Unfortunately, there $x$ has to be away from the reflection hyperplanes associated with all roots for some results, and
some results hold for particular root systems only and are too weak.
\end{remark}

On the other hand, the first limit in \eqref{limit-1f1} can be used to derive the following central limit theorem for $N=1$:

\begin{theorem}\label{clt-dunkl-drift}
  For a $1$-dimensional Dunkl process  $(X_t^\lambda)_{t\ge0}$ with  $k>0$,
  drift $\lambda\ne0$,  and start in 0, the process $(X_t^\lambda-t\lambda)/\sqrt t$ tends in distribution
to the normal distribution $N(0,1)$ for $t\to\infty$.
\end{theorem}

\begin{proof}  By symmetry we may assume $\lambda>0$. Using  \eqref{density-transition-dunkl-drift-gen} and
  $c_k=\frac{1}{2^{k+1/2}\Gamma(k+1/2)}$, we see that
$(X_t^\lambda-t\lambda)/\sqrt t$ has the density
 \begin{equation}
 \phi_t(y):= \frac{1}{2^{k+1/2}\Gamma(k+1/2)}  \frac{e^{-\lambda^2 t/2}}{t^{k}}  e^{-(y\sqrt t+\lambda t)^2/(2t)} E_k(y\sqrt t+\lambda t,\lambda)
\cdot |y\sqrt t+\lambda t|^{2k}.
 \end{equation}
Hence, by \eqref{limit-1f1} and the doubling formula for the gamma function, 
$$ \lim_{t\to\infty} \phi_t(y)= \frac{\Gamma(2k+1)}{2^{2k+1/2}\Gamma(k+1/2)\Gamma(k+1)}e^{-y^2/2}=\frac{1}{\sqrt{2\pi}}e^{-y^2/2}$$
locally uniformly in $y\in\mathbb R$. This means that  the probability densities $\phi_t$ tend vaguely and thus weakly to $N(0,1)$ as claimed.
\end{proof}

\begin{remark} By Theorem 2 of \cite{RdJ} there exists a limit result for $E_k(ty,\lambda)$ for$N\ge2$,
  $t\to\infty$ and $y,\lambda$ in the interior
  of the same Weyl chamber $C\subset \mathbb R^N$, which corresponds to the first formula in \eqref{limit-1f1}.
This result corresponds to a short-time asymptotics for the Dunkl heat kernels; see Corollary 2  of \cite{RdJ}.
We conjecture that the following variant of Theorem 2 of \cite{RdJ} holds:
\begin{equation}\label{conjecture-dunkl-kernel-limit}
  \lim_{t\to\infty} E_k(y\sqrt t+\lambda t, \lambda)
  \cdot (2\pi)^{N/2} c_k e^{-\langle y\sqrt t+\lambda t, \lambda\rangle} \sqrt{w_k( y\sqrt t+\lambda t)w_k(\lambda)}=1
\end{equation}
for all $\lambda,y\in\mathbb R^N$, where $\lambda$ is  contained in the interior of some Weyl chamber.

Conjecture \eqref{conjecture-dunkl-kernel-limit} for some $\lambda$ and the proof of the CLT above  immediately imply that 
 for a Dunkl process  $(X_t^\lambda)_{t\ge0}$ with  this
  drift $\lambda\ne0$  and start in 0, $(X_t^\lambda-t\lambda)/\sqrt t$ tends in distribution
  to the $N$-dimensional standard normal distribution $N(0,I_N)$.

  We point out that the CLT with a normal distribution as limit  usually does not hold for $\lambda$ on the boundary of  some Weyl chamber.
  As a counterexample take the Weyl group $\{\pm1\}^2$ on $\mathbb R^2$ and $\lambda=(1,0)$. In this case,  the CLT above holds for the first
  component of $(X_t^\lambda)$ while in the second component no centering occurs, and the limit distribution of the second component
  is a two-sided Gamma-distribution.
\end{remark}

\section{Limits for Bessel processes with drift}

In this section we  study moments,  moment functions,  SLLNs, and CLTs for
Bessel processes and  hybrid Dunkl-Bessel processes with drift with a focus on Bessel processes.
In  both cases, the computations in the proof of Lemma \ref{expectation-dunkl-drift-xi}
do not lead to meaningful results (except for $k=0$). For this reason we immediately turn to analogues of 
the moment functions in Definition \ref{def-moments-dunkl} which are adapted to both types of processes:

\begin{definition}\label{def-moments-Bessel}
Define the
  multivariate (modified)
  moment functions of order 1 and 2 by
  $$m_1^{Bessel}(x):=(m_{1;1}^{Bessel}(x),\ldots,m_{1;N}^{Bessel}(x)):= \frac{1}{J_k(x,\lambda)}\nabla_\lambda J_k(x,\lambda)$$
  and
  $$m_2^{Bessel}(x):=(m_{2;j,l}^{Bessel}(x))_{j,l=1,\ldots,N}:= \frac{1}{J_k(x,\lambda)}H_\lambda J_k(x,\lambda)$$
  for $x\in\mathbb R^N$ with the gradient $\nabla_\lambda$ and the Hessian $H_\lambda$ w.r.t.
  $\lambda\in\mathbb R^N$.
\end{definition}

With this definition, all assertions of Lemma  \ref{bounds-moments-dunkl}, Propositions \ref{mart-m1}, \ref{mart-m2}, and \ref{SLLN-m1},
 and Lemma  \ref{limit-m1-Dunkl}
 remain valid with the same proofs for Bessel processes with drift and the moment functions in Definition \ref{def-moments-Bessel}.
In particular, for $N=1$ and Bessel processes with drift $\lambda>0$, the strong law \ref{SLLN-1dim-dunkl}
 and the CLT \ref{clt-dunkl-drift} remain valid with the same proofs (in fact, some arguments are easier here, as no jumps appear):

 \begin{theorem}\label{SLLN-CLT-1dim-bessel}
  For a $1$-dimensional Bessel process  $(X_t^\lambda)_{t\ge0}$ on $[0,\infty[$ with continuous paths, $k\ge 0$,
    drift $\lambda>0$,  and start in 0,
    $$\lim_{t\to\infty} \frac{X_t^\lambda}{t} = \lambda \quad{a.s.},$$
  and  $(X_t^\lambda-t\lambda)/\sqrt t$ tends in distribution
to the normal distribution $N(0,1)$.
 \end{theorem}

 These results follow also from the limit theorems for Sturm-Liouville hypergroups on $[0,\infty[$ with exponential growth in
     Section 7.3 of \cite{BH}; see also \cite{Z1, Z2} for the original results.

On the other hand, the strong law \ref{SLLN-1dim-dunkl}
 and the CLT \ref{clt-dunkl-drift} are not valid for hybrid Dunkl-Bessel processes with start in $0$, as here 
 the distributions are invariant under the action of the Weyl group. For this reason there are no particular limit theorems except for those
 which appear from the fact that the
 projections of these hybrid processes
 onto some Weyl chamber $C$ under $p_C$ just are corresponding Bessel processes. For this reason we now study Bessel processes only.

\begin{example}\label{example-11}
  Let us collect all data   for the simplest nontrivial example with $N=1$, 
  $k=1$ and,  w.l.o.g., $\lambda=1$. This example clarifies all assertions above for $N=1$.
  In this case, by Lemma \ref{int-rep-Dunkl-1dim-Prp} and a short computation, 
$$c_1=(2\pi)^{-1/2},  \quad\quad   E_1(\lambda,x)= \frac{e^{\lambda x}}{\lambda x} -\frac{ \sinh(\lambda x)}{\lambda^2 x^2} ,
\quad\quad J_1(\lambda,x)=\frac{\sinh(\lambda x)}{\lambda x},$$
$$m_1^{Dunkl}(x)= \frac{(x^2-x)e^x -x \cosh x +2 \sinh x}{xe^x- \sinh x},\quad\quad  m_1^{Bessel}(x)= \frac{x \cosh x -\sinh x}{\sinh x}, $$
$$ L^{1, Dunkl}_1 f(x)= \frac{1}{2}f^{\prime\prime}(x)  +\frac{f^{\prime}(x)}{x} \frac{x^2e^x -x \cosh x + \sinh x}{xe^x- \sinh x}
+ \frac{f(-x)-f(x)}{x^2} \cdot \frac{-xe^{-x}+ \sinh x}{xe^x- \sinh x},$$
$$ L^{1, Bessel}_1 f(x)= \frac{1}{2}f^{\prime\prime}(x)  +\cth x f^{\prime}(x),$$
and
$$ \tilde L^{1, Dunkl}_1 f(x)= \frac{1}{2}f^{\prime\prime}(x)  +\cth x f^{\prime}(x)
+ \frac{f(-x)-f(x)}{2x^2} $$
for the generators $ L^{1, Bessel}_1$ and $\tilde L^{1, Dunkl}_1 $ of the Bessel and hybrid Dunkl-Bessel processes respectively.
Furthermore,  if $(X_t^1)_{t\ge0}$ is a Dunkl process $(X_t^1)_{t\ge0}$  with drift 1 and start in $0$,
 then, by \eqref{density-transition-dunkl-drift-gen}, $X_t$ (for $t>0$) has the Lebesgue density
 $$\frac{ e^{- t/2}}{t^{3/2} \sqrt{2\pi} } e^{-y^2/(2t)} (ye^y-\sinh y) \quad\quad (y\in\mathbb R).$$
 Notice that  these densities also lead to the CLT \ref{clt-dunkl-drift}  and
 to the CLT-part in Theorem \ref{SLLN-CLT-1dim-bessel}.
\end{example}

The limit results in Theorem \ref{SLLN-CLT-1dim-bessel} for $N=1$ are closely related to corresponding
limit results for the geometric examples discussed in Subsections 3.1-3.4, i.e., for particular multiplicities $k$ depending of the root system.
In fact, in these geometric cases, Bessel processes with and without drift on some Weyl chamber $C\subset\mathbb R^N$
with start in $0$ can be regarded as projections of Brownian motions with and without drift on  associated 
Euclidean vector spaces $(V, \langle\,.\,,\,.\,\rangle)$ of finite dimension as in Subsection 3.1. We here mainly have the following
examples in mind for the fields $\mathbb F=\mathbb R,\mathbb C, \mathbb H$ with real dimension $d=1,2,4$:
\begin{enumerate}\itemsep=-1pt
  \item[\rm{(1)}] Type $A_{N-1}$: $k=d/2=1/2,1,2$,  $V=\mathbb H_N(\mathbb F)$;
 \item[\rm{(2)}]  Type $B_{N}$: $k_1=(M-N+1)d/2-1/2, \> k_2=d/2$, $V=M_{M,N}(\mathbb F)$ ($M\ge N$ an integer).
\end{enumerate}
Now let $g:V\to C$ be the canonical projection, and $\lambda\in C$ a fixed drift. Let $\tilde\lambda\in V$ with $g(\tilde\lambda)=\lambda$.
The classical SLLN for Brownian motions $(B_t^{\tilde \lambda})_{t\ge0}$ on $V$ with drift $\tilde\lambda$ and start in $0\in V$
 implies the SLLN for the associated Bessel processes
$(X_t^\lambda)_{t\ge0}$ on $C$ with start in 0 and with drift $\lambda$, i.e.,
 $X_t/t\to\lambda$ for $t\to\infty$  a.s..
 
Moreover, for all $t>0$, $ \sqrt t(B_t^{\tilde \lambda}/t-\tilde\lambda)=(B_t^{\tilde \lambda}-t\tilde\lambda)/\sqrt t$ is $N(0,I)$-distributed on
$V$ with the $dim \> V$-dimensional identity matrix $I$.
We now assume that $\lambda$ is in the interior of $C$. Then $g$ is continuously differentiable in $\tilde\lambda$, and the
Jacobi matrix $J_g(\tilde\lambda)$ satisfies  $J_g(\tilde\lambda)J_g(\tilde\lambda)^T=I_N$ with $I_N$ the identity matrix of dimension $N$;
see for instance Theorem 5.17 in Ch.~1 of \cite{He} and the text before where it is stated that
$J_g(\tilde\lambda)=(I_N, 0)\in \mathbb R^{N\times dim \> V}$ in a hidden way. For the $A_{N-1}$-case an elementary proof of these statements can be found e.g.~in Section 5.3 of \cite{D}.
We thus conclude from the Delta-method for  CLTs of transformed random variables  (see Section  3.1 of \cite{vV})
that $X_t^{\lambda}-t\lambda)/\sqrt t$ converges in distribution to $N(0,I_N)$. In summary:

\begin{theorem}\label{SLLN-CLT-geometric-bessel}
  Let $R$ be a root system and $k\ge0$ a multiplicity such that the associated Dunkl-Bessel theory is adapted  to some geometric case.
  Let $(X_t^\lambda)_{t\ge0}$ be a Bessel process on some associated closed Weyl chamber $C\subset \mathbb R^N$
  with some drift $\lambda\in C$ with start in 0.
Then
    $$\lim_{t\to\infty} \frac{X_t^\lambda}{t} = \lambda \quad{a.s.}.$$
Moreover, if  $\lambda$ is in
  the interior of $C$, then
   $(X_t^\lambda-t\lambda)/\sqrt t$ tends in distribution
to  $N(0,I_N)$.
 \end{theorem}

The CLT does not hold for  $\lambda$ on the boundary of $C$ for obvious geometric reasons.

We also notice that in  the geometric cases $A_{N-1}$, $k=1$, our Bessel processes with drift
$\lambda=(-N+1, -N+3,\ldots, N-3, N-1)\in \mathbb R^N$ can be also regarded as projections of Brownian motions without drift on
the symmetric spaces $SL(N,\mathbb C)/SU(N,\mathbb C)$ (see \cite{AuV, RV3} and references there) which means that 
these examples may be seen as examples of  Heckman-Opdam processes associated with the root systems $A_{N-1}$ for which
there are SLLNs and CLTs in \cite{V2}. In a similar way, there are geometric cases of type $B_{N}$ where certain Bessel processes
with drift may be seen as examples of  Heckman-Opdam processes associated with the root systems  $BC_{N}$.

\begin{remark} We conjecture that the SLLN and CLTs in the settings of Theorems \ref{SLLN-CLT-1dim-bessel} and \ref{SLLN-CLT-geometric-bessel}
   hold for general root systems and 
   multiplicities $k\ge0$ for drift vectors $\lambda$ in the interior of the corresponding Weyl chambers. In fact, it is shown in \cite{V3} that at least the CLT
   holds in general under mild conditions on  $\lambda$ and the starting point.

Moreover, the SLLN and the CLT can possibly be derived
  from limit results for the  Bessel functions which correspond to the limit \eqref{limit-1f1}  for one-dimensional Dunkl kernels
  which was shown  above by using classical results for  $_1F_1$-functions from \cite{NIST}.
  This limit result can be also derived directly from the integral representation for Dunkl kernels in
  Proposition \ref{int-rep-Dunkl-1dim-Prp}. We expect that further  explicit integral representations for Bessel functions
  for $N\ge2$ lead to analogues of \eqref{limit-1f1} and thus  to the SLLN and CLT in  Theorem
  \ref{SLLN-CLT-geometric-bessel}.
  In fact, for the case $A_{N-1}$ ($N\ge2$) and $k>0$ there exists a recursive  integral representation in \cite{Am}
  based on  corresponding formulas on Jack polynomials in \cite{OO}. Moreover,  \cite{R4} contains an
  integral representation for the root systems $B_N$ with $k=(k_1,k_2)=(k_1, d/2)$ with $d=1,2,4$ and $k_1>0$ sufficiently large.
For $k_1=(M-N+1)\cdot d/2-1/2$ with $M\ge N$ an integer, this case is  geometric  where 
the SLLN and CLT hold by  \ref{SLLN-CLT-geometric-bessel}; see  \cite{R4} and the discussion in Section \ref{geometric-b}.
\end{remark}

\begin{remark}
The spectra of square roots of Wishart processes on the cones of positive definite $N\times N$-matrices over
$\mathbb F=\mathbb R, \mathbb C, \mathbb H$ (see \cite{Bru, DDMY})
are Bessel processes of type $B_N$ without drift with  $k=(k_1,k_2)=(k_1, d/2)$ with $d=1,2,4$ and $k_1$ real and sufficiently large.
In the literature there exist several concepts of  Wishart processes with drifts; namely the concept studied in \cite{DDMY}
which generalizes the case of one-dimensional Bessel processes with drift in \cite{PY}, as well as another different concept
in  \cite{GMM} and references cited there.
The latter concept of drifts does not seem to be adapted to our drift concept for Bessel processes while the first one does, i.e.,
spectra of  Wishart processes with drifts as in \cite{DDMY} are Bessel processes with drift of type B.
  \end{remark}

\begin{remark} Limit theorems for Bessel processes $(X_t^{(k_1,k_2)})_{t\ge0}$ of type $B_N$ for fixed $t>0$, $k_2\ge0$
  and $k_1\to\infty$ can be found in \cite{RiV}. \end{remark}

\section{Declaration}

This version was cerrected by using ChatGPT5.6 mainly in view of misprints and notational errors.
Parts of the proof of Theorem 3.4(1) were corrected and
rewritten by using ChatGPT5.6.

\end{document}